\documentclass[11pt,a4paper,reqno]{amsart}
\usepackage{amsmath,amssymb,graphics,epsfig,color,enumerate}

\usepackage[T1]{fontenc}
\usepackage[utf8]{inputenc}

\usepackage[dvipsnames]{xcolor}
\usepackage{dsfont}  
\usepackage[normalem]{ulem}
\usepackage{verbatim}
\usepackage{xcolor}
\usepackage{hyperref}
\definecolor{refkey}{gray}{.75}
\definecolor{labelkey}{gray}{.5}

\newtheorem{Theorem}{Theorem}[section]
\newtheorem{Fact}{Fact}

\newtheorem{Lemma}[Theorem]{Lemma}
\newtheorem{Proposition}[Theorem]{Proposition}
\newtheorem{Corollary}[Theorem]{Corollary}
\newtheorem{Remark}[Theorem]{Remark}

\newtheorem{Definition}[Theorem]{Definition}

 \definecolor{darkgreen}{rgb}{0,0.4,0}

\definecolor{light}{gray}{0.9}

 \let\so=\o

\newcommand{\cA}{\ensuremath{\mathcal A}}

\newcommand{\cD}{\ensuremath{\mathcal D}}
\newcommand{\cE}{\ensuremath{\mathcal E}}
\newcommand{\cF}{\ensuremath{\mathcal F}}
\newcommand{\cG}{\ensuremath{\mathcal G}}
\newcommand{\cH}{\ensuremath{\mathcal H}}

\newcommand{\cL}{\ensuremath{\mathcal L}}

\newcommand{\cN}{\ensuremath{\mathcal N}}

\newcommand{\cV}{\ensuremath{\mathcal V}}

\newcommand{\bbC}{{\ensuremath{\mathbb C}} }

\newcommand{\bbE}{{\ensuremath{\mathbb E}} }

\newcommand{\bbG}{{\ensuremath{\mathbb G}} }

\newcommand{\bbI}{{\ensuremath{\mathbb I}} }

\newcommand{\bbL}{{\ensuremath{\mathbb L}} }

\newcommand{\bbN}{{\ensuremath{\mathbb N}} }

\newcommand{\bbP}{{\ensuremath{\mathbb P}} }

\newcommand{\bbR}{{\ensuremath{\mathbb R}} }

\newcommand{\bbT}{{\ensuremath{\mathbb T}} }

\newcommand{\bbX}{{\ensuremath{\mathbb X}} }

\newcommand{\bbZ}{{\ensuremath{\mathbb Z}} }

    \let\d=\delta  \let\e=\varepsilon
 \let\g=\gamma       \let\l=\lambda
      \let\o=\omega      
  \let\s=\sigma \let\t=\tau   
 \let\x=\xi \let\z=\zeta
     \let\L=\Lambda 
\let\O=\Omega      

\newcommand{\rosso}{\textcolor{black}}

\newcommand{\be}{\begin{equation}}
\newcommand{\en}{\end{equation}}

\newcommand{\ra}{\rangle}
\newcommand{\la}{\langle}

\newcommand{\mfm}{\mathfrak{m}}

\newcommand{\da}{\downarrow}

\newcommand{\rmd}{{\rm d}}

\newcommand{\dive}{{\rm div}\hspace{0.04cm}}

\author[A.~Faggionato]{Alessandra Faggionato}
\address{Alessandra Faggionato.
  Department of Mathematics, Sapienza University of Rome. 
  P.le Aldo Moro 2, 00185 Rome, Italy}
\email{faggiona@mat.uniroma1.it}

\author[M.~Salvi]{Michele Salvi}
\address{Michele Salvi.
  University of Rome Tor Vergata, Via della ricerca scientifica 1, 00133,  Rome. Italy}
\email{salvi@mat.uniroma2.it}

\title[Scaling limit of the complex mobility matrix]{Scaling limit of the complex mobility matrix for  the  random conductance model on $\bbT^d_N$}

\begin{document}

\begin{abstract}

We consider a continuous-time random walk on the $d$-dimensional torus $\bbT^d_{N}=\bbZ^d/N \bbZ^d$, possibly with long-range, but finite, jumps. The law of the jumps is regulated by a random environment $\xi$ yielding a stationary and ergodic field of random conductances. The complex mobility matrix  $\sigma_N^\xi(\omega)$  measures the linear response of the random walk to a $\cos(\o t)$--type oscillating external field. By investigating the homogenization properties of the medium, and assuming in addition that the conductances have finite second moment, we show that, for almost every realization of the environment $\xi$, the complex mobility matrix $\sigma_N^\xi(\omega)$ converges as $N\to+\infty$ to a deterministic limiting matrix $\sigma(\omega)$ and provide different characterizations of   $\sigma(\omega)$.

\smallskip

   
\smallskip


\end{abstract}

\maketitle


\bigskip

\section{Introduction}

The complex mobility matrix describes the linear transport response of a medium to a weak, time--periodic external field. It is also known as the AC conductivity matrix -- where AC stands for alternating current -- and is a fundamental concept in solid-state physics and materials science  (cf.~\cite{AM,KTH}). 
More precisely,  consider an external field
oscillating in time as $\cos (\omega t)$,  where $\o>0$ denotes  the field frequency. If the field direction is given by the unit vector $\hat v$, then for a typical ohmnic system the macroscopic instantaneous current varies in time as   \[\Re \left( {\rm e}^{i \o t}  \sigma (\o) \hat v \right)\l + o(\l)\,,
\]
where $\l$ is the  field intensity  and $\s(\o)$ is the complex mobility matrix. Above $\Re(\cdot)$ denotes the real part of a complex number.

\smallskip

For disordered media, $\s(\o)$ is expected to satisfy  universal laws  both in the low-frequency 
and high-frequency regimes, possibly at low temperatures  (cf.~e.g. \cite{ABSO}, \cite{Dyre}, \cite[Chapter~6]{POF},\cite{Sah}). To study these laws, one considers very large systems in order to neglect boundary effects, which mathematically corresponds to taking the infinite-volume limit.
The linear response behavior has been rigorously studied and derived in  \cite{FM,FMS,FSi,JPS}    for several stochastic processes, including cases with zero spectral gap, for which a perturbative approach is precluded. In this work, we provide a first rigorous contribution to the analysis of the infinite-volume limit of the complex mobility matrix in disordered media, focusing on the specific case of a random walk with random conductances. The latter, also known as the conductance model, is a paradigmatic model in probability theory for studying motion in disordered media.

\smallskip

To perform the infinite-volume limit, we assume that a random conductance field on $\bbZ^d$,
with finite-range conductances, is assigned in a stationary and ergodic manner. Given a scaling parameter $N\in \{1,2,\dots\}$, we then consider the 
periodized conductance field with unit cell $[0,N)^d\cap \bbZ^d$.
  This naturally defines a random conductance field on the discrete torus $\bbT^d_N:=\bbZ^d/N \bbZ^d$, and the corresponding continuous--time random walk on $\bbT^d_N$, for which the conductances coincide with the jump probability rates.
In \cite{FSi} the authors have provided a characterization (recalled in Fact~\ref{ale_vitt} and Fact \ref{ale_vitt2} below) of the complex mobility matrix $\s_N^\xi (\o)$  of this random walk on $\bbT^d_N$ under   a $\cos(\o t)$-like time--dependent drift,  where $\xi$ denotes a realization of the original  conductance field on $\bbZ^d$ or - more generally -  the random \emph{environment}.

\smallskip
  
Starting from the  formula provided in \cite{FSi}, and under the additional assumption of finite second moment of the conductances, we prove that, for almost all realizations $\xi$, the complex mobility matrix $\s_N^\xi (\o)$ converges as $N\to +\infty$ to a 
deterministic matrix $\s(\o)$, for which we provide an intrinsic  characterization. This is the content of our limit theorems (cf.~Theorem~\ref{teo1} and Theorem \ref{teo1bis}).  In Section~\ref{sec_alg_id} we also provide an alternative characterization of the limiting matrix $\s(\o)$, as well as  of the real and imaginary parts of $a \cdot \s(\o) a$ for $a\in \bbR^d$.

\smallskip

According to our   formulas for  $\s(\o)$, this matrix formally corresponds to the complex mobility matrix of the Markov process \emph{environment viewed from the particle} when perturbed by a weak $\cos (\omega t)$-like oscillating external field.
We stress that this correspondence is only formal, since a full linear response analysis of this Markov process is still lacking. This is a non-trivial task, as the process has zero spectral gap --      precluding perturbative arguments -- while approaches based on regeneration times, as used in the derivation of the Einstein relation for random walks or diffusions in random environments
 (cf.~\cite{GGN,GMP,MP}),  do not apply because the perturbed system can have zero time-averaged macroscopic velocity.

 We also note that, in general,  the complex mobility matrix introduced above for 
 $\o>0$ reduces to the  standard mobility matrix
  when $\o=0$. By the Einstein relation this coincides with   the diffusion matrix. 
 We recall that  the  quenched infinite volume limit of the diffusion matrix for the random walk in random environment, as the size $N$ of the periodization diverges, has been derived in \cite{CI} and \cite{O} under much more restrictive conditions on the conductance field, due to the absence of the complex--massive term $i\o$ present in our settings (as  will become clear from the proof). The same quenched limit, but for  diffusions in random environments, was derived in \cite{BP}.

\smallskip

We conclude by providing some comments on our derivation. We have avoided any uniform ellipticity assumption for the conductances and we have used the method of $2$--scale convergence to be able to cover a large range of random conductance fields. This method  
was  introduced by Nguetseng  in homogenization theory of partial differential equations in \cite{Nu} and further developed by Allaire in \cite{A}. It has   been successfully adapted to stochastic homogenization,   allowing to deal also with  random differential operators on singular structures (cf. \cite{ZP} and references therein) and then with random conductance models on the supercritical percolation cluster \cite{MPrw} or more general random graphs on $\bbR^d$ (cf.~\cite{F_hom} and, for resistor networks, \cite{F_resistor}). 
The 2-scale convergence method guarantees  some weak compactness  of $L^2$--bounded functions
where weak convergence is understood as the convergence of the integral of these functions against observables expressed in terms of both the macroscopic position of the random walk (when embedding  $\bbT^d_N$ into $\bbT^d:=\bbR^d/\bbZ^d$) 
and the microscopic behavior of the environment around the walker. Note the presence of two scales: the macroscopic scale and the microscopic scale. This compactness is formalized in Lemma~\ref{compatto1} and Lemma \ref{compatto2} below and relies on the ergodicity of the disordered medium. As a consequence, for typical environments and at the cost of extracting a subsequence $(N_k)$ of the scaling parameter $N$,   we can show that the functions  appearing in the representation of $\s_N^\xi(\o)$ provided by \cite{FSi}, as well as their discrete gradients, 
have a limit (cf.~Corollary~\ref{caciotta}). These functions correspond to the functions $\theta^\xi_N$ introduced in Section~\ref{sec_dim_I} which satisfy the discrete elliptic system \eqref{olleboiccic}.
On the other hand, this  limit could depend on $\xi$ itself (even  the subsequence $(N_k)$ could  depend on $\xi$). The main effort is then to prove the uniqueness of the limit, at least for the averaged part, 
which effectively determines
 $\s_N^\xi(\o)$ (cf.~Lemma~\ref{grovis} and Corollary~\eqref{gioioso}). 
 Our proof requires us to assume the conductances to have finite second moment, but we point out that many intermediate results are also valid for conductances with finite mean only (see Remark \ref{remarkone} for a list of places where the finite second moment condition is used \rosso{in Section~\ref{sec_dim_I}}).

 In the  discussion above, for the sake of clarity, we have  neglected the effects of boundary conditions. Since the complex mobility matrix
 $\s_N^\xi(\o)$  refers to the random walk on the discrete torus, one must indeed consider periodic boundary conditions, which impose certain global constraints. To investigate the stochastic homogenization of the medium, in Section \ref{sec_dim_I}  we first analyze the system in the ``bulk'', avoiding the effect of the boundaries, while in Section \ref{sec_conclusione}    we combine the bulk stochastic homogenization with an analysis of the periodic boundary conditions. Finally, we point out that our derivation was
 inspired by the  approach     developed in \cite{F_hom} and  \cite{F_resistor}, from which we  benefited in obtaining  some preliminary results and  techniques, partly  collected in  Section~\ref{sec_hom}.

Our limit theorem is  of qualitative nature. Quantitative results, in the  spirit of recent developments of stochastic homogenization (cf. e.g. \cite{AKM,GNO} and references therein), will be obtained under more restrictive assumptions and will be presented separately.

\medskip

{\bf Outline of the paper}. In Section \ref{sec_mod_res} we introduce the model, we recall the characterization  for $\s_N^\xi(\o)$ provided in \cite{FSi} (cf.~Fact~\ref{ale_vitt} and Fact~\ref{ale_vitt2})   and present our main results. We treat first the case of nearest neighbor jumps and then, in Section \ref{rotterdam}, we consider the more general case of long-range jumps. The main result  is given by Theorem \ref{teo1} in the nearest neighbor case and by the stronger Theorem \ref{teo1bis} in the general case. 

\smallskip

In Section \ref{sec_hom} we collect some preliminary definitions and results, some of which are already known (cf.~\cite{F_hom,F_resistor}). In particular, in Section \ref{sec_typ_env} we  define the full-measure set $\Omega_{\rm typ}$ of well-behaving environments to which we will restrict for all our results.

\smallskip

Section \ref{sec_dim_I} is the core of the proof of Theorem \ref{teo1bis}.  We analyze in the bulk the microscopic version, called $\theta_\varepsilon^\zeta$ and $\nabla_\varepsilon\psi_\varepsilon^\zeta$, of some quantities appearing in the formula for the complex mobility matrix $\sigma_N^\xi(\o)$ appearing in Fact \ref{ale_vitt2}. We show that they are uniformly bounded in $L^2$ (Lemma \ref{lemma_stime}) and deduce that a.s. they 2-scale converge along a subsequence to some non-trivial limit as $\varepsilon$ goes to $0$ (Corollary \ref{caciotta}). In Lemma \ref{fvg}, Lemma \ref{venezia} and Lemma \ref{biauzzo}  we collect equations satisfied by the limiting objects.

\smallskip

In  Section \ref{sec_conclusione} we conclude the proof of Theorem \ref{teo1bis}. This is easily achieved thanks to the results of Section \ref{sec_dim_I} after one takes care of the boundary terms induced by the torus geometry for finite $N$ (Lemma \ref{pirano}).

\smallskip

Finally, in Appendix~\ref{proof_prop_zap} we prove Proposition~\ref{zap} by means of Dirichlet forms theory. This proposition is used in Section  \ref{sec_mod_res}  to  define the limiting matrix  $\s(\o)$ and also in Section \ref{sec_dim_I} to get some uniqueness of the limit points in the 2-scale convergence  after some partial averaging (see Corollary~\ref{gioioso}).

\section{Model and main result}\label{sec_mod_res}

For the sake of clarity, we first introduce the model and state the main results in the simpler case of nearest-neighbour jumps. In Section \ref{rotterdam} we will extend our analysis to random walks that can perform long-range jumps.
So, for the moment, let  $\cN:=\{z\in\bbZ^d:\,|z|=1\}$ be the set of neighbours of the origin. We denote by $\bbE_d$ the non-oriented edges of $\bbZ^d$ and write $e_1,\dots, e_d$ for the canonical basis of $\mathbb R^d$.

\smallskip

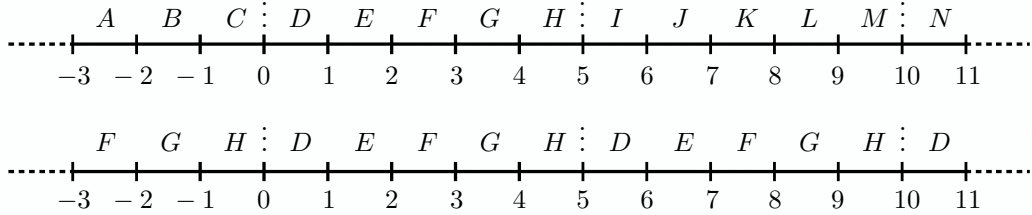
\begin{figure}[h!]
	\begin{center}
		\setlength{\unitlength}{0.12\textwidth}  
		\begin{picture}(7.5,2.1)(0,0)
		\put(-0.4,0){\includegraphics[scale=0.35]{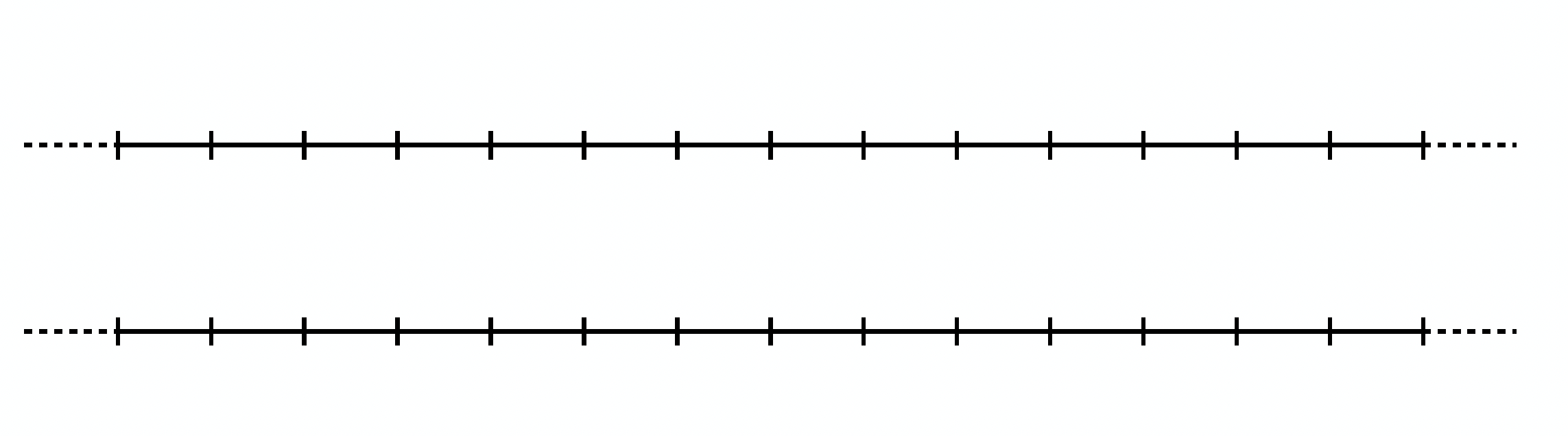}}
		\put(0,1.25){\small \, $-3\; \  -2\; \  -1 \hspace{0.55cm} 0\hspace{0.68cm} 1\hspace{0.65cm} 2\hspace{0.68cm} 3\hspace{0.65cm} 4\hspace{0.66cm} 5\hspace{0.67cm} 6\hspace{0.68cm} 7\hspace{0.65cm} 8\hspace{0.65cm} 9\hspace{0.65cm} 10\hspace{0.45cm} 11$}
		\put(0,0.27){\small \, $-3\; \  -2\; \  -1 \hspace{0.55cm} 0\hspace{0.68cm} 1\hspace{0.65cm} 2\hspace{0.68cm} 3\hspace{0.65cm} 4\hspace{0.66cm} 5\hspace{0.67cm} 6\hspace{0.68cm} 7\hspace{0.65cm} 8\hspace{0.65cm} 9\hspace{0.65cm} 10\hspace{0.45cm} 11$}
		\put(0.2,1.7){\small $\quad A\hspace{0.6cm}B\hspace{0.57cm}C\hspace{0.55cm}D\hspace{0.55cm}E\hspace{0.55cm}F\hspace{0.55cm}G\hspace{0.55cm}H\hspace{0.57cm}I\hspace{0.63cm}J\hspace{0.6cm}K\hspace{0.54cm}L\hspace{0.55 cm}M\hspace{0.52cm}N$}
		\put(1.7,1.7){\color{black}\vdots}
		\put(4.19,1.7){\color{black}\vdots}
		\put(6.7,1.7){\color{black}\vdots}
		\put(0.2,0.72){\small $\quad {\color{black}F\hspace{0.58cm}G\hspace{0.55cm}H}\hspace{0.55cm}D\hspace{0.55cm}E\hspace{0.55cm}F\hspace{0.55cm}G\hspace{0.55cm}H\hspace{0.55cm}{\color{black}D\hspace{0.55cm}E\hspace{0.55cm}F\hspace{0.55cm}G\hspace{0.55cm}H\hspace{0.55cm}D\hspace{0.55cm}}$}
		\put(1.7,0.72){\color{black}\vdots}
		\put(4.19,0.72){\color{black}\vdots}
		\put(6.7,0.72){\color{black}\vdots}
		\end{picture}
	\end{center}
\caption{Periodization procedure in 1d  for nearest-neighbour conductances. On the first line:  conductances 
			$(c_{x,x+1})_{x\in\mathbb Z}$ (e.g.~$c_{0,1}=D$). On the second line: periodized conductances  $(c^{(5)}_{x,x+1})_{x\in\mathbb Z}$.
	} \label{figure111} 
	\end{figure}

\begin{figure}
\includegraphics[scale=0.2]{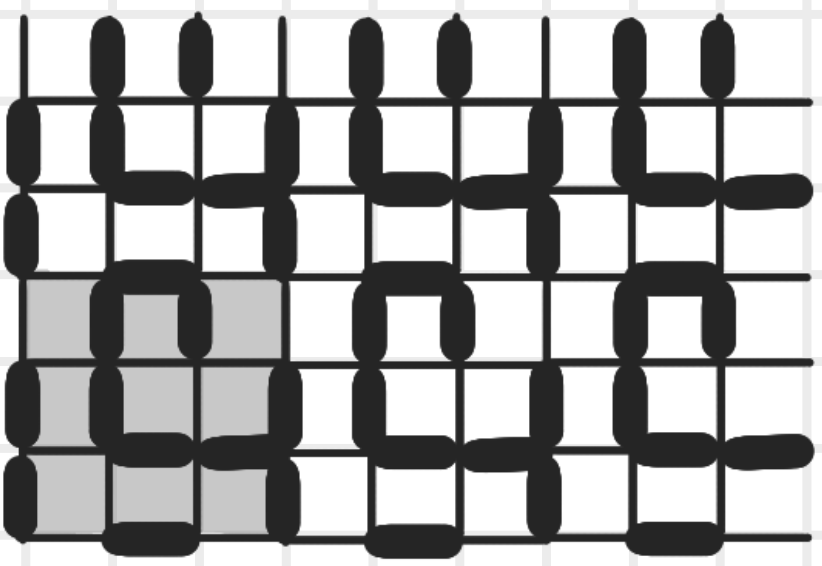}
\caption{Portion of a periodized environment in $\bbZ^2$ with $N=3$. The grey region corresponds to $[0,N)^d$. Different bond thicknesses represent different conductances.}
\label{figure3}
\end{figure}

Consider  a probability space $(\O,\cF, \bbP)$.  We write $\bbE[\cdot]$ for the expectation with respect to $\bbP$. Elements $\xi\in \O$ are called \emph{environments}.
For $\xi\in\Omega$, we associate to each pair of neighbouring points $x,x+z\in \bbZ^d$, with $z\in\cN$,  a  positive symmetric jump rate, or \emph{conductance}, $c_{x,x+z}(\xi)$. Symmetry means here that $c_{x,x+z}(\xi)= c_{x+z, x} (\xi)$.
For a given $N\in \bbN_+:=\{1,2,\dots\}$ and a realization of conductances $\{c_{x,x+z}(\xi)\}_{x\in\bbZ^d,z\in\cN}$, we can define its $N$-periodized version as the set of conductances $\{c_{x,x+z}^{(N)}(\xi)\}_{x\in\bbZ^d,z\in \cN}$ 
such that $c_{x,x+e_j}^{(N)}(\xi)=c_{x,x+e_j}(\xi)$ for all $x\in   [0,N)^d \cap \bbZ^d$ and $j\in\{1,\dots,d\}$ and $c_{x,x+e_j}^{(N)}(\xi)=c_{x+y,x+y+e_j}^{(N)}(\xi)$ for all  $y\in N \bbZ^d$ (see Figure \ref{figure111} and Figure \ref{figure3}).

We introduce  the canonical projection 
\[ 
\pi_N:\bbZ^d \to \bbT^d_N,
\]
where  $\bbT^d_N:=\bbZ^d/N \bbZ^d$ is the discrete torus. As usual, given $x\in \bbT^d_N$ and $z\in \bbZ^d$, we write $x+z$ for the element $x+\pi_N(z) \in \bbT^d_N$. Moreover, we extend the definition of conductances to pairs of neighbouring vertices $x,y\in \bbT^d_N$ in the torus by setting   $c_{x,y}^{(N)}(\xi):=c^{(N)}_{x',y'}(\xi)$, where $x',y'$ is any pair of neighbouring points in $\bbZ^d$ such that $\pi_N(x')=x$ and $\pi_N(y')=y$.

\smallskip

Following \cite{FSi} we now introduce the complex mobility matrix for the random walk on $\bbT^d_N$ with jump rates given by the random conductances  (also called \emph{random conductance model}, cf.~\cite{Bi}). Notice that the discussion and results in \cite{FSi} are for a generic \textit{fixed} family of   conductances, so in this part  the randomness of the conductances will play no role.

Given $\xi \in \O$ and $N\in \bbN_+$, 
we consider the continuous-time random walk $(X^{\xi,N}_t)_{t\geq 0}$ on $\bbT^d_{N}=\bbZ^d/N \bbZ^d$
with nearest-neighbour jumps such that, for $x\in\bbT^d_{N}$ and $z\in\cN$, 
\be\label{unpert_rw}
{\rm P}( X^{\xi,N}_{t+\rmd t} =x+z\,| \,X^{\xi,N}_t=x)= c^{(N)} _{x,x+z}(\xi) {\rm d}t+o({\rm d}t) \,.
\en
Since the conductances are strictly positive, the random walk is irreducible and the
 uniform distribution $\mfm_N$ on $\bbT^d_{N}$ is reversible for this process.
We now perturb the random walk $(X^{\xi,N}_t)_{t\geq 0}$ by putting a  time--oscillatory external field in a given direction. To this aim we fix 
\begin{itemize}
	\item $v\in \bbR^d\setminus\{0\}$,  the direction\footnote{Strictly speaking, the term direction refers to vectors $v$ with $|v|=1$. Since normalization of $v$
 is not required for our results, we use the same terminology for all  $v\in \bbR^d\setminus\{0\}$.} of the field;
	\item $\o>0$, the frequency of the field;
	\item $\l>0$, the intensity of the field.
\end{itemize} 
We consider the 
continuous-time random walk $(X^{\xi,N,\l}_t)_{t\geq 0}$ on $\bbT^d_{N}$
with nearest-neighbour jumps such that, for $x\in\bbT^d_N$ and $z\in\cN$,  
\be\label{pert_rw}
{\rm P}( X^{\xi,N,\l}_{t+{\rm d}t} =x+z\,|\, X^{\xi,N,\l}_t=x)
	= c^{(N)} _{x,x+z}(\xi) \exp \{ \l \cos(\o t ) z \cdot v\}  {\rm d} t+o({\rm d}t) \,.\en
	We also call $T:= 2\pi /\o$ the period of the field. 
By irreducibility  and aperiodicity, the discrete-time Markov chain $\big(X^{\xi,N,\l}_{nT}\big)_{n\geq 0}$ admits a unique invariant distribution $\mfm_{\xi,N,\l}$ on $\bbT^d_{N}$. It is simple to check that  $\mfm_{\xi,N,\l}$ is the unique initial distribution such that the law on the Skorohod space of c\`adl\`ag paths of the perturbed random walk is left invariant by time-translations which are integer multiples of the period $T $. This  law is called \emph{oscillatory steady state} (OSS) of the perturbed random walk.

By the Dynkin martingale, the mean instantaneous velocity at time $t$ of the perturbed random walk in the OSS is given by
\begin{equation}\label{mean_velocity}
V^{\xi,N,\l}(t)
= \sum_{x\in \bbT^d_{N}}  \sum_{z\in \cN} {{\text{P}}^\text{OSS}} ( X^{\xi,N,\l}_t=x)\,c^{(N)} _{x,x+z}(\xi) \exp \{ \l \cos(\o t ) z \cdot v\} z\,.
\end{equation}
The probability $ {{\text{P}}^\text{OSS}}$ refers to the OSS, i.e.~to the perturbed random walk starting with distribution $\mfm_{\xi,N,\l}$. 

\smallskip 

We fix once and for all some notation frequently used in the rest of the paper. Given a probability space $(\bbX,\mu)$ we denote by $L^2(\bbX,\mu)$ or $L^2(\mu)$ the Hilbert space with scalar product $\la f,g\ra_{L^2(\mu)} := \int \overline{f (x)} g(x) \rmd \mu(x)$. We denote by $L^2_R(\bbX,\mu)$ or $L^2_R(\mu)$ its closed subspace of real functions. In general, for $L^p$ spaces we will use the subindex $R$ to indicate that we deal with real $L^p$--functions.
Given an operator $A: L^2_R(\mu)\ni f\mapsto Af \in L^2_R(\mu)$ defined over  real functions,  we can extend $A$  to complex functions by considering the operator with domain $\cD(A)+ i \cD(A) \subset L^2(\mu)$ mapping $f+i g$ into $Af+i Ag$ for all $f,g\in \cD(A)\subset L^2_R(\mu)$. This  new operator on complex functions will still be denoted by $A$, but we will specify that it is the complex extension.

\smallskip
Let us recall some linear response results from  \cite{FSi}. We denote  by  $\cL^{\xi}_N: L^2(\mfm_N)\to L^2(\mfm_N) $  the extension to complex functions of the Markov generator of the unperturbed random walk in $L^2_R(\mfm_N)$:
\be\label{cancellami}
\cL^{\xi}_N f(x) 
	:= \sum_{z\in\cN } c^{(N)}_{x,x+z}(\xi) \big( f(x+z)-f(x) \big)\qquad \forall x \in \bbT^d_{N}\,.
\en
Note that $ \cL^{\xi}_N$  
is a symmetric operator on $L^2(\mfm_N)$, hence it can be diagonalized along an orthonormal basis.  By  irreducibility of the unperturbed random walk, zero is a simple eigenvalue of  $\cL^{\xi}_N$, while all other eigenvalues are strictly negative. This implies that the operator $i\o - \cL^{\xi}_N$ is invertible and $(i\o - \cL^{\xi}_N)^{-1} f= \int_0^{+\infty} {\rm e}^{ -(i \o - \cL^{\xi}_N) s} f{\rm d} s $  for any function $f$ with 
zero mean $\mfm_N(f)$
(note that the integral is absolutely convergent in $L^2(\mfm_N)$ since the integrand has exponentially decaying  norm in  $L^2(\mfm_N)$). 

Recall that the frequency $\o $ is positive  and that the direction of the field  $v\in \bbR^d$ appears in the rates of the perturbed random walk \eqref{pert_rw}. Below we denote by  $\Re (a)$ and $ \Im(a)$ the real and imaginary part of $a\in\bbC$, respectively. The following fact is stated and proved in  \cite{FSi}:

\begin{Fact}\cite[Theorem~5.1, Remark~5.2]{FSi} \label{ale_vitt}
Consider the random walk $(X^{\xi,N,\l}_t)_{t\geq 0}$ on $\bbT^d_{N}$ with nearest-neighbour jumps.	There is a unique  $d\times d$ complex  matrix $\s^{\xi}_{N}(\o )$, called   {\bf complex mobility matrix},  such that 
 \begin{equation}\label{stilton}
\partial_{\l=0} V^{\xi,N,\l} (t)   
	= \Re \left( {\rm e}^{i \o t} \s^{\xi}_{N}(\o ) v \right)  \qquad  \forall v \in \bbR^d,\; \; t\geq 0\,.
\end{equation}
Let  $c^{\xi}_N$ and  $\g^{\xi}_N$ be the  functions from $ \bbT^d_{N} $ to $ \bbR^d$ given by  
\begin{align}
 c^{\xi}_N(x)
 	&:=\sum_{j=1}^d \bigl[ c^{(N)}_{x,x+e_j}(\xi)+ c^{(N)}_{x,x-e_j}(\xi)\bigr] e_j\,, \label{def_ciccio_N1}\\
\gamma^{\xi}_N(x) 
	&:=\sum_{j=1}^d  \bigl[ c^{(N)}_{x,x+e_j}(\xi)- c^{(N)}_{x,x-e_j}(\xi)\bigr] e_j
	= \sum_{z\in\cN} c^{(N)}_{x,x+z}(\xi) z \label{def_gamma_N1}
	\,.
\end{align}
Then, $\s^{\xi}_N(\o)$ is the symmetric matrix such that for $j,k\in \{1,2,\dots, d\}$ it holds 
 \begin{equation}\label{jabba2}
\begin{split}
\s^{\xi}_N(\o)_{j,k}
	& =  \mfm_N [c^{\xi}_{N,j}]\d_{j,k}-2 \la \gamma_{N,j}^{\xi}, (i \o -\cL^{\xi}_N)^{-1} \gamma _{N,k}^{\xi} \ra_{L^2(\mfm_N)} \\
	&=  \mfm_N [c_{N,j}^{\xi }]\d_{j,k}-2 \int_0^{+\infty}  \la \gamma_{N,j}^{\xi}, {\rm e}^{ -(i \o - \cL^{\xi}_N) s} \gamma_{N,k}^{\xi}\ra_{L^2(\mfm_N)} \;{\rm d}s
	\end{split}
\end{equation}
where $c^{\xi}_{N,j}$ and $\gamma_{N,j}^{\xi}$ denote the $j$--th coordinate of $c^\xi_N$ and $\g^\xi_N$, respectively.

\end{Fact}

\begin{Remark}\label{rem_segno}The  law of the perturbed random walk  does not change if we substitute $\o$ by $-\o$, see \eqref{pert_rw}. In particular, by extending the above treatment also to $\o<0$, by the invariance of the mean instantaneous velocity we would get  $\Re \big( {\rm e}^{i \o t} \s^{\xi}_{N}(\o ) v \big) =\Re \big( {\rm e}^{-i \o t} \s^{\xi}_{N}(-\o ) v \big) $ for all $t\geq 0$ and $v$, thus implying that  $\s^\xi(-\o)=\overline{\s^\xi(\o)}$.
\end{Remark}

We focus our interest on the infinite volume limit of the complex mobility matrix when $\xi$ is random. We assume that the group $\bbZ^d$ acts on $\Omega$  and we  write  $(\t _g)_{g\in \bbZ^d}$ for this action. We recall that this means that $\t_g:\O\to\O$ is measurable for any $g\in\bbZ^d$, $\t_0={\rm id}$ and $\t_g\circ \t_{g'}=\t_{g+g'}$ for all $g,g'\in\bbZ^d$.
We require the probability measure $\bbP$ on the environment space $\O$ to be stationary and ergodic with respect to this action. Stationarity means that  $\bbP(A)=\bbP(\t_g A)$ for all $A\in\cF$ and $g\in\bbZ^d$, while ergodicity means that  $\bbP(A)\in\{0,1\}$ for any translation invariant set $A\in\cF$  (i.e.~such that $\t_g A=A$ for all $g\in\mathbb Z^d$).  We additionally require that $L^2(\bbP)$ is separable (equivalently, that $L^2_R(\bbP)$ is separable)  and for  all $x,y,z\in \bbZ^d$ we impose the covariant relation
\be\label{covariante}
c_{x+z,y+z}(\xi)= c_{x,y}(\t_z\xi)\,.
\en
A typical choice  is given by  $\O=(0,+\infty) ^{\bbE_d}$, $\cF$ being the Borel $\s$-algebra. The action maps are given by  $\t_z \xi:= (\xi_{e-z}: e\in \bbE_d)$ for $\xi=(\xi_e:e\in\bbE_d)$ and $e-z:=\{x-z,y-z\}$ for $e=\{x,y\}$. Moreover conductances are given by $c_{x,y}(\xi)=\xi_{\{x,y\}}$ (usually, one just writes $\xi_{x,y}$). Identity \eqref{covariante} is then trivially satisfied. One has to require stationarity and ergodicity of $\bbP$ (e.g.~$\bbP$ could be a product probability measure).
\smallskip

Before stating our main theorem, we need a technical result in order to introduce the infinite-volume counterpart $\bbL$ of the Markov generator appearing in \eqref{cancellami} and some of its properties.
We recall (cf.~\cite[Section~1.1]{FOT}) that a \emph{symmetric form}  $\cE$ on $L^2_R(\bbP)$ is  a positive semi-definite symmetric bilinear form  defined on a dense subspace $\cD(\cE)$  of $L^2_R(\bbP)$. It is called \emph{closed} if $\cD(\cE)$ is complete w.r.t.~the norm $\|\cdot \|_1$ defined  on $\cD(\cE)$ as $\|f\|_1^2:= \|f\|_{L^2(\bbP)}^2+ \cE(f,f)$ (i.e. if $\|\cdot\|_1$--Cauchy sequences converge). Finally  we recall that $\cE$ is a Dirichlet form if it is a Markovian   closed symmetric form. By the discussion in \cite[Section~1.1]{FOT}, the Markovian property for a closed symmetric form is equivalent to the following: if $u\in \cD(\cE)$ and $v$ is a normal contraction of $u$, i.e.
\be\label{zompicchia}
|v(\xi)-v(\xi')| \leq |u(\xi)-u(\xi')|\;\; \forall \xi,\xi'\in \O, \qquad |v(\xi) |\leq |u(\xi) | \;\; \forall \xi \in \O\,,
\en
then $v\in \cD(\cE)$ and $\cE(v,v)\leq \cE(u,u)$.

\begin{Proposition}\label{zap}   Suppose that the probability measure $\bbP$ on $\O$ is stationary and ergodic w.r.t.~the action of the group $\bbZ^d$ and that the covariant relation \eqref{covariante} is satisfied. In addition suppose  that $\bbE[c_{0,z}(\xi)]<+\infty$ for all $z\in\cN$.
Consider the bilinear form 
\be \label{def_dir_form}
\cE(f,g):= 
	\frac{1}{2} \sum _{z\in\cN} 
	\bbE\big[ c_{0,z}(\xi) \left( f (\t_z \xi) - f(\xi)\right)\, \left( g (\t_z \xi) - g(\xi)\right) \big]\,,
\en
with domain $ \cD(\cE):=\big\{ f\in L^2_R(\bbP)\,:\, \sum_{z\in\cN}\bbE\big[ c_{0,z}(\xi) \left( f (\t_z \xi) - f(\xi)\right)^2\big]<+\infty\big\}$. Then  $\cE$ is a  Dirichlet form and measurable bounded functions  form a core of $\cE$, i.e. a dense  subset of $\cD(\cE)$ w.r.t. $\|\cdot\|_1$.

Moreover,  there exists a unique negative semidefinite self-adjoint operator $\bbL: \cD(\bbL)\to L^2_R(\bbP)$ with $\cD(\bbL)\subset L^2_R(\bbP)$  such that 
\be\label{istria}
 \cD(\cE)= \cD(\sqrt{-\bbL})\,, \qquad \| \sqrt{-\bbL}   f \| ^2_{L^2(\bbP)} = \cE(f,f) \qquad \forall f\in \cD(\cE)\,.
\en
 For  any function $f\in L^2_R(\bbP)$ such that
\begin{itemize}
\item[(i)]  the map $\xi\mapsto c_{0,z}(\xi) \left( f(\t_z \xi) - f(\xi) \right)$ belongs to $L^1_R(\bbP)$ for any $z\in \cN$,
\item[(ii)]  $L f(\xi):=\sum_{z\in \cN} c_{0,z}(\xi) \left( f(\t_z \xi) - f(\xi) \right)$ belongs to $L^2_R(\bbP)$, \end{itemize}
it holds $f\in\cD(\bbL)$ and  $\bbL f=Lf $.
 In particular, if $\bbE[c_{0,z}^2]<+\infty$ for all $z\in\cN$, then   any  bounded and measurable  function $f:\O\to\bbE$  belongs to $\cD(\bbL)$ and it holds $\bbL f=L f$. 
\end{Proposition}

We postpone the proof of Proposition~\ref{zap} to  Appendix~\ref{proof_prop_zap}.

We can finally state our main result in the nearest-neighbour setting, concerning the infinite volume limit of the complex mobility matrix. The theorem, in which we recap the assumptions presented above, relies on the homogenization properties of the medium.
\begin{Theorem}\label{teo1}
Suppose that the probability measure $\bbP$ on $\O$ is stationary and ergodic w.r.t.~the action of the group $\bbZ^d$ and that $L^2(\bbP)$ is separable.  In addition, suppose that the conductance field satisfies   the covariant relation \eqref{covariante}  and  that $\bbE[c^2_{0, z}]<+\infty$ for all $z\in\cN$. Then, for $\bbP$--a.a.~$\xi \in \O$ it holds 
\be\label{santa_pizza}
\lim_{N\to +\infty} \s^{\xi}_N(\o)= \s (\o)\,,
\en 
where the symmetric $d\times d$ complex matrix $\s (\o)$ is defined as follows.   Let  $c$ and  $\g$ be the  functions from $ \O $ to $ \bbR^d$ given by  \begin{align}
& c(\xi)
	:=\sum_{j=1}^d \bigl[ c_{0,e_j}(\xi)+ c_{0, -e_j}(\xi)\bigr] e_j\,, \label{def_ciccio}\\
& \gamma(\xi) 
	:=\sum_{j=1}^d \bigl[ c_{0,e_j}(\xi)- c_{0, -e_j}(\xi)\bigr] e_j
	= \sum_{z\in\cN} c_{0,z}(\xi) z
\,.\label{def_gamma}
\end{align}
Then  $\s(\o)$ is the symmetric matrix with entries
 \begin{equation}\label{jabba2bis}
\begin{split}
\s (\o)_{j,k}& =  \bbE [c_j]\d_{j,k}-2 \la \gamma_j, (i \o -\bbL)^{-1} \gamma_k \ra_{L^2(\bbP) }\\
&=  \bbE [c_j ]\d_{j,k}-2 \int_0^{+\infty}\la \gamma_j , {\rm e}^{ -(i \o - \bbL) s} \gamma_k \ra_{L^2(\bbP)}\, {\rm d}s \,,
\end{split}
\end{equation}
for $j,k\in \{1,2,\dots, d\}$\,, where $\bbL$ is the self-adjoint operator defined in Proposition \ref{zap} (extended to complex functions).
\end{Theorem}

\begin{Remark}\label{nevischio} By spectral calculus, see~\cite[Lemma~2.5]{FM},  the 
 operator 
$i \o -\bbL$ is invertible with bounded inverse operator $(i \o -\bbL)^{-1} $ and  for any $f\in L^2(\bbP)$ with zero mean it holds  $ (i \o -\bbL)^{-1} f =\lim_{S\to+\infty}\int_0^{S}{\rm e}^{ -(i \o - \bbL) s} f \rmd s$ where the limit is in $L^2(\bbP)$. We point out that  $\bbE[\g_k]=0$.  In particular, the second identity in \eqref{jabba2bis} is just an application of \cite[Lemma~2.5]{FM}.
\end{Remark}

 
 \subsection{Extension to the long-range model}\label{rotterdam}
Theorem \ref{teo1} can be extended to the case of long-range, but finite, interactions. The model is essentially the same as in the previous section, but now the set $\cN$ of neighbours, i.e.~the set of jumps the random walk is allowed to take, is more general. To be precise, we define the following:

\begin{Definition}[Good neighbourhood] \label{def:neighborhood}
	We say that a set $\cN\subset \bbZ^d\setminus\{0\}$ is a good neighbourhood of the origin if it satisfies the following properties:
	\begin{itemize}
		\item $\cN$ is finite;
		\item $\cN$ is symmetric (that is, if $z\in \cN$  then $-z\in \cN$);
		\item for all $x\not=y$ in $\bbZ^d$ there are $x_0:=x, x_1, \dots, x_{k-1}, x_k:=y$ with  $x_{i+1}-x_i\in \cN$. 
	\end{itemize}
\end{Definition}

We choose once and for all a good neighbourhood $\cN$ of the origin. For any point $x\in\bbZ^d$ we call the set $x+\cN$ the neighbourhood of $x$. Moreover, we fix once and for all a subset $\cN_*\subset \cN$ such that $\cN$ can be partitioned as $\cN=  \cN_* \sqcup(- \cN_*)$. 
For example, the previous construction  in the nearest-neighbour case can be recovered by  choosing  $\cN_*$ as  the canonical basis $\{e_1,\dots,e_d\}$.

\smallskip

Consider now a probability space $(\O,\cF, \bbP)$. 
For a given $\xi\in\Omega$, we associate to each pair of points $x,y \in \bbZ^d$ a random non-negative symmetric jump rate, or conductance, $c_{x,y}(\xi)$ such that  $c_{x,y}(\xi)>0$ if and only if $y-x\in \cN$. 

\smallskip

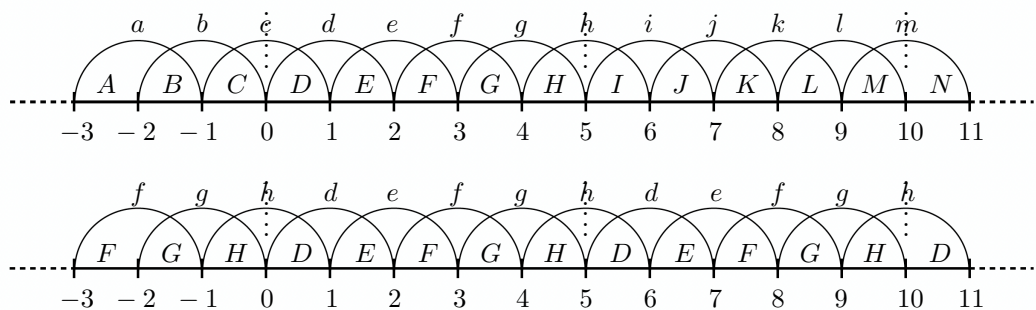
\begin{figure}[h]
	\begin{center}
		\setlength{\unitlength}{0.12\textwidth}  
		\begin{picture}(7.5,2.1)(0,0)
		\put(-0.4,0){\includegraphics[scale=0.318]{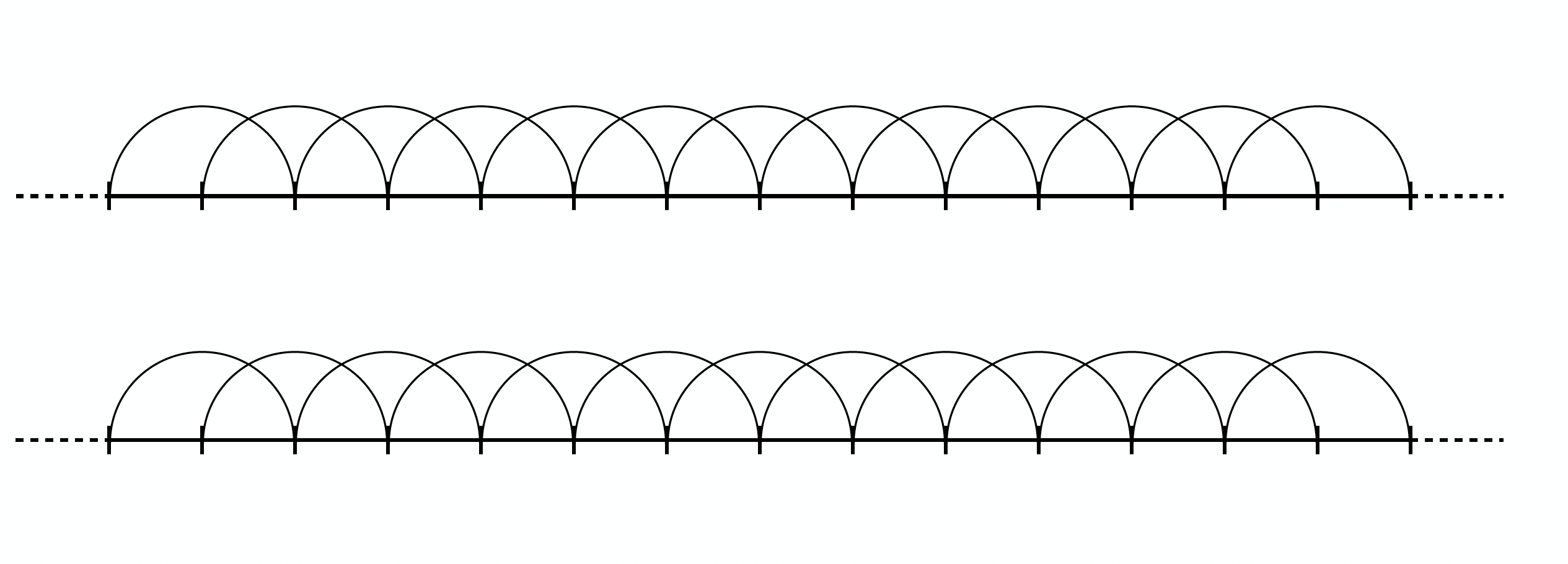}}
		\put(-0.03,1.68){\small \, $-3\; \  -2\; \  -1 \hspace{0.55cm} 0\hspace{0.68cm} 1\hspace{0.65cm} 2\hspace{0.68cm} 3\hspace{0.65cm} 4\hspace{0.66cm} 5\hspace{0.67cm} 6\hspace{0.68cm} 7\hspace{0.65cm} 8\hspace{0.65cm} 9\hspace{0.65cm} 10\hspace{0.45cm} 11$}
		\put(-0.03,0.36){\small \, $-3\; \  -2\; \  -1 \hspace{0.55cm} 0\hspace{0.68cm} 1\hspace{0.65cm} 2\hspace{0.68cm} 3\hspace{0.65cm} 4\hspace{0.66cm} 5\hspace{0.67cm} 6\hspace{0.68cm} 7\hspace{0.65cm} 8\hspace{0.65cm} 9\hspace{0.65cm} 10\hspace{0.45cm} 11$}
		\put(0.15,2.05){\small $\quad A\hspace{0.6cm}B\hspace{0.57cm}C\hspace{0.55cm}D\hspace{0.55cm}E\hspace{0.55cm}F\hspace{0.55cm}G\hspace{0.55cm}H\hspace{0.57cm}I\hspace{0.63cm}J\hspace{0.6cm}K\hspace{0.54cm}L\hspace{0.55 cm}M\hspace{0.52cm}N$}
		\put(1.66,2.5){\color{black}\vdots}
		\put(1.66,2.25){\color{black}\vdots}
		\put(4.16,2.5){\color{black}\vdots}
		\put(4.16,2.25){\color{black}\vdots}
		\put(6.67,2.5){\color{black}\vdots}
		\put(6.67,2.25){\color{black}\vdots}	
		\put(0.15,0.72){\small $\quad {\color{black}F\hspace{0.58cm}G\hspace{0.55cm}H}\hspace{0.55cm}D\hspace{0.55cm}E\hspace{0.55cm}F\hspace{0.55cm}G\hspace{0.55cm}H\hspace{0.55cm}{\color{black}D\hspace{0.55cm}E\hspace{0.55cm}F\hspace{0.55cm}G\hspace{0.55cm}H\hspace{0.55cm}D\hspace{0.55cm}}$}
		\put(1.66,1.17){\color{black}\vdots}
		\put(1.66,0.92){\color{black}\vdots}
		\put(4.16,1.17){\color{black}\vdots}
		\put(4.16,0.92){\color{black}\vdots}	
		\put(6.67,1.17){\color{black}\vdots}
		\put(6.67,0.92){\color{black}\vdots}	
		\put(0.42,2.53){\small $\quad a\hspace{0.68cm}b\hspace{0.68cm}c\hspace{0.66cm}d\hspace{0.67cm}e\hspace{0.66cm}f\hspace{0.65cm}g\hspace{0.68cm}h\hspace{0.65cm}i\hspace{0.7cm}j\hspace{0.68cm}k\hspace{0.66cm}l\hspace{0.67 cm}m\hspace{0.52cm}$}
		\put(0.419,1.2){\small $\quad {\color{black}f\hspace{0.64cm}g\hspace{0.67cm}h}\hspace{0.64cm}d\hspace{0.64cm}e\hspace{0.66cm}f\hspace{0.65cm}g\hspace{0.68cm}h\hspace{0.65cm}{\color{black}d\hspace{0.68cm}e\hspace{0.63cm}f\hspace{0.65cm}g\hspace{0.67 cm}h\hspace{0.52cm}}$}
		\end{picture}
	\end{center}
\caption{Periodization procedure in 1d   for nearest-neighbours and second-nearest-neighbours conductances ($\cN=\{\pm e_1, \pm 2 e_1\}$ and  $\mathcal N_*=\{e_1, 2e_1\}$). On the first line:  conductances 
				$(c_{x,x+w})_{x\in\mathbb Z,\,w\in \cN_*}$ (e.g.~$c_{0,1}=D$ and $c_{0,2}=d$). On the second line:   periodized conductances  $(c^{(5)}_{x,x+w})_{x\in\mathbb Z,\,w\in\mathcal N_*}$.  
		}
	\label{figure2} 
\end{figure}

We can now define the periodized rates (see Figure~\ref{figure2}).
\begin{Definition}[Periodized rates] \label{def_per_rates}Given $\xi \in \O$ and $N\in \bbN_+$  we define the periodized rates $\big\{c_{x,y}^{(N)}(\xi)\big\}_{x,y\in\bbZ^d}$ as follows:
\begin{itemize}
\item if $y-x\not \in\cN$, then $c_{x,y}^{(N)}(\xi):=0 $,
\item if $y-x\in\cN$, then $c_{x,y}^{(N)}(\xi)$ equals  the rate  $c_{u,w}$ where the pair $u,w$ is the unique one such that  $u\in [0,N)^d\cap \bbZ^d$, $w-u\in \cN_*$ and $\{u,w\}= \{x,y\} + N z$ for some $z\in \bbZ^d$. 
\end{itemize}
\end{Definition}

Equivalently, the second item above could be replaced by the following: for any $x\in \bbZ^d$ and $z\in\cN_*$ we set $c_{x,x+z}^{(N)}(\xi):=c_{u,u+z}(\xi)$ where $u$ is the unique element in $[0,N)^d\cap\bbZ^d$ such that $x\in u+N \bbZ^d$.
We point out the following properties: 
$ c_{x,y}^{(N)} (\xi) = c_{y,x} ^{(N)} (\xi)$ for all $ x,y\in\bbZ^d$, and 
$ c_{x,y}^{(N)}(\xi)= c_{x+ N z , y+ Nz} ^{(N)} (\xi)$ for all $x,y,z \in \bbZ^d$.

\vspace{1.5cm}

The modelization follows then exactly as in the previous section
 when
  \be\label{rose}
 N>2\|\cN\|_\infty\,,  \qquad \|\cN\|_\infty:= \max\{\|z\|_\infty: z\in \cN\}\,.
  \en
 Indeed, in this case one can extend the definition of conductances on the discrete torus $\bbT^d_N$ as follows. Recall the canonical projection $\pi_N:\bbZ^d\to \bbT^d_N$.  Take $x\not=y$ in $\bbT^d_N$.  If $x=\pi_N(x')$, $y=\pi_N(y')$ and $y'=x'+w$ for some $x',y'\in \bbZ^d$ and  $w\in \cN$, then $w$ is univocally determined as $N>2 \|\cN\|_\infty$. In this case we set $c_{x,y}^{(N)}(\xi):=c_{x',y'}^{(N)}(\xi)$ (note that the r.h.s.~does not depend on the particular pair $x',y'$). In all other cases we set $c_{x,y}^{(N)}(\xi):=0$. Then one introduces the unperturbed random walk $(X^{\xi,N}_t)_{t\geq 0}$ on $\bbT^d_N$ with jump rates given by equation   \eqref{unpert_rw} with $z\in \cN$. This random walk 
 is irreducible (due to the last property listed in Definition~\ref{def:neighborhood}) and  the uniform measure $\mfm_N$   on $\bbT^d_N$ is a  reversible distribution. Finally,  one introduces the perturbed random walk $(X^{\xi,N,\l}_t)_{t\geq 0}$ on $\bbT^d_N$  with jump rates given by equation \eqref{pert_rw} with $z\in \cN$. The definition of   the oscillatory steady state (OSS) is the same as in the nearest-neighbour case.
The instantaneous velocity at time $t$ in the OSS is again given by \eqref{mean_velocity} and the
complex extension of the Markov generator of the unperturbed random walk by \eqref{cancellami}, where in both cases we consider of course the new neighbouring set $\cN$. Fact \ref{ale_vitt} can be replaced by the following:
\begin{Fact}\cite[Theorem~5.3]{FSi} \label{ale_vitt2}
	Consider the random walk $(X^{\xi,N,\l}_t)_{t\geq 0}$ on $\bbT^d_{N}$ that can take jumps in a good neighbourhood set $\cN$ as in Definition \ref{def:neighborhood}. Then there exists a unique  $d\times d$ complex  matrix $\s^{\xi}_{N}(\o )$, called   {\bf complex mobility matrix}, satisfying \eqref{stilton}.  $\s^{\xi}_{N}(\o )$ can be characterized as follows.
	Let   $\g^{\xi}_N$ be the  function from $ \bbT^d_{N} $ to $ \bbR^d$ given by  
	\begin{align}
\gamma^{\xi}_N(x) 
&:=\sum_{z\in\cN_*}  \bigl[ c^{(N)}_{x,x+z}(\xi)- c^{(N)}_{x,x-z}(\xi)\bigr] z
= \sum_{z\in\cN} c^{(N)}_{x,x+z}(\xi) \,z \label{def_gamma_N}
\,.
\end{align} 
	Then for $j,k\in \{1,2,\dots, d\}$ it holds 
	\begin{equation}\label{jabba2*}
	\begin{split}
	\s^{\xi}_N(\o)_{j,k}
	& =  \sum_{z\in\mathcal N} \mfm_N[c_{\cdot,\cdot+z}^{(N)}(\xi)]z_jz_k 
	-2 \la \gamma_{N,j}^{\xi}, (i \o -\cL^{\xi}_N)^{-1} \gamma _{N,k}^{\xi} \ra_{L^2(\mfm_N)} \\
	&= \sum_{z\in\mathcal N} \mfm_N[c_{\cdot,\cdot+z}^{(N)}(\xi)]z_jz_k 
	-2 \int_0^{+\infty}  \la \gamma_{N,j}^{\xi}, {\rm e}^{ -(i \o - \cL^{\xi}_N) s} \gamma_{N,k}^{\xi}\ra_{L^2(\mfm_N)} \;{\rm d}s
	\end{split}
	\end{equation}
	where 
	$\gamma_{N,j}^{\xi}$ denotes the $j$--th coordinate of 
	$\g^\xi_N$.
\end{Fact}

Remark~\ref{rem_segno} remains valid in this extended context.
Proposition \ref{zap}   holds verbatim whenever  $\cN$ is a good neighbourhood as in Definition \ref{def:neighborhood} (its  proof in Appendix~\ref{proof_prop_zap} is indeed given in the general case). Finally we restate the main theorem in the long-range case.

\begin{Theorem}\label{teo1bis} Suppose that the probability measure $\bbP$ on $\O$ is stationary and ergodic w.r.t.~the action $(\t_z)_{z\in \bbZ^d}$ of the group $\bbZ^d$ and that $L^2(\bbP)$ is separable.  In addition, suppose that the conductance field satisfies   the covariant relation \eqref{covariante}  and   that $\bbE[c^2_{0, z}]<+\infty$ for all $z\in\cN$, with $\cN$ a good neighbourhood as in Definition \ref{def:neighborhood}. Then for $\bbP$--a.a.~$\xi \in \O$ it holds 
	\be\label{santa_pizza_bis}
	\lim_{N\to +\infty} \s^{\xi}_N(\o)= \s (\o)\,,
	\en 
	where the symmetric $d\times d$ complex matrix $\s (\o)$ is defined as follows.   Let  
	$\g$ be the  function from $ \O $ to $ \bbR^d$ given by  \begin{align}
	 \gamma(\xi) 
	:=\sum_{z\in\cN_{*}} \bigl[ c_{0,z}(\xi)- c_{0, -z}(\xi)\bigr] z
	= \sum_{z\in\cN} c_{0,z}(\xi) z
	\,.\label{def_gammaV}
	\end{align}
	Then  $\s(\o)$ is the symmetric matrix with entries
\begin{equation}\label{jabba2bisV}
	\begin{split}
	\s (\o)_{j,k}
		& =  \sum_{z\in\mathcal N} \bbE[c_{0,z}]z_jz_k-2 \la \gamma_j, (i \o -\bbL)^{-1} \gamma_k \ra_{L^2(\bbP) }\\
		&= \sum_{z\in\mathcal N} \bbE[c_{0,z}]z_jz_k-2 \int_0^{+\infty}\la \gamma_j , {\rm e}^{ -(i \o - \bbL) s} \gamma_k \ra_{L^2(\bbP)}\, {\rm d}s \,,
	\end{split}
	\end{equation}
for $j,k\in \{1,2,\dots, d\}$, where $\bbL$ is the self-adjoint operator defined in Proposition \ref{zap} (extended to complex functions) with the present good neighbourhood $\cN$.

\end{Theorem}
Remark \ref{nevischio} holds verbatim  referred to \eqref{jabba2bisV} instead of \eqref{jabba2bis}.

\smallskip

Theorem~\ref{teo1bis} is a special fluctuation-dissipation theorem, since it provides a representation of  the linear response of the system  to a weak external field in terms of equilibrium quantities. Trivially Theorem~\ref{teo1bis} implies Theorem~\ref{teo1}. We devote  the next sections to its proof.  Note that   the   external field no longer appears in Theorem~\ref{teo1bis}, as is usual in fluctuation-dissipation theorems. As a consequence, it will not appear in the next sections either.  In particular, we warn the reader that  in what follows  the letter $v$ (used above to denote the field direction)  will be used to denote other objects. 

\subsection{Some algebraic identities satisfied by $\s(\o)$}\label{sec_alg_id}
In this section we provide some alternative expressions for the limiting complex mobility matrix $\sigma(\o)$, valid both in the nearest-neighbor case with $\cN_*$ given by the canonical basis and in the more general case treated in Section~\ref{rotterdam}. Given a  function $u:\O\to \bbC$  we define the gradient of $u$ as the  function $\nabla u: \O \times \cN \to \bbC$ such that 
 \begin{equation}\label{cantone}
 \nabla u (\xi, z):= u (\t_z \xi)-u (\xi)\,.
 \end{equation}
Furthermore, we define the function  $\theta:\Omega\to\bbC^d$ with entries given by 
$$\theta_k=(i\o-\bbL)^{-1}\gamma_k$$
for each $k=1,\dots,d$. The above equation has to be thought in $L^2(\bbP)$ (recall the invertibility of the operator $i\o-\bbL$ stated in Remark~\ref{nevischio}).

\begin{Proposition}
An alternative expression for the matrix $\sigma(\o)$ in Theorem \ref{teo1bis} is 
\be \label{acqua3}
    \s (\o)_{j,k}
        = 2 \sum_{z\in\mathcal N_*}\bbE \big[ c_{0,z} z_j\big( z_k + \nabla \theta _k(\cdot,z)\big) \big]\,.
\en
Equivalently, in matrix form, 
$$
\sigma(\o)
    =2\sum_{z\in\mathcal N_*}\bbE\big[\Lambda(z)\big(\mathcal Z(z)+\Psi(z)\big)\big]\,,
$$
where $\Lambda_{jk}(z):=c_{0,z}(\xi)z_j\delta_{jk}$, $\mathcal Z_{jk}(z):= z_k$ and $\Psi_{jk}(z) :=\nabla\theta_k(\xi,z)$.

\smallskip

These formulas further simplify in the nearest-neighbours case of Theorem \ref{teo1}, where they become
$$
\s (\o)_{j,k}= 2 \bbE \big[ c_{0,e_j}\big( \delta_{jk} + \nabla \theta _k(\cdot,e_j)\big) ]
$$
and 
\begin{equation}\label{labour}
\sigma(\o)
    =2\bbE\big[\Lambda\big(\bbI +\Psi\big)\big]\,,
\end{equation}
with $\Lambda_{jk}=c_{0,e_j}\delta_{jk}$ and $\Psi_{jk}=\nabla\theta_k(\xi,e_j)$. 
\end{Proposition}
\begin{proof}
   We restrict our attention to \eqref{acqua3}, since the remaining formulas are an easy consequences of it. Due to  \eqref{jabba2bisV} we just need to show that 
   \begin{align}
       & \sum_{z\in\mathcal N} \bbE[c_{0,z}]z_jz_k=2 \sum_{z\in\mathcal N_*}\bbE \big[ c_{0,z}] z_j z_k\,,\label{pallina1}\\
       & \la \gamma_j, (i \o -\bbL)^{-1} \gamma_k \ra_{L^2(\bbP)}=-\sum_{z\in\mathcal N_*}\bbE \big[ c_{0,z} z_j \nabla \theta _k(\cdot,z)\big]\,. \label{pallina2}
   \end{align}
   Eq.~\eqref{pallina1} follows easily since $\cN=\cN_*\sqcup (-\cN_*)$ and $\bbE[c_{0,z}]=\bbE[c_{-z,0}]=\bbE[c_{0,-z}]$  by the covariance relation \eqref{covariante} implying that  $c_{-z, 0} (\t_z \xi)=c_{0,z}(\xi)$, the stationarity of $\bbP$ and  the symmetry of the conductances. By definition of $\theta_k$, equation \eqref{pallina2} reads   $\la \gamma_j, \theta_k \ra_{L^2(\bbP)}=
   \sum_{z\in \cN}\bbE[c_{0,z} z_j \theta_k]_{L^2(\bbP)}=-\sum_{z\in\mathcal N_*}\bbE \big[ c_{0,z} z_j \nabla \theta _k(\cdot,z)\big]$.  This follows immediately by noticing that
    $\bbE[c_{0,-z}z_j\theta_k]=\bbE[c_{0,z}z_j\theta_k(\tau_z\cdot)]$ (to check this last identity use that  $c_{0,-z}(\t_z \xi)\theta_k(\t_z \xi)= c_{0,z}( \xi)\theta_k(\t_z \xi)$ and argue as above).
\end{proof}
We point out that \eqref{labour} corresponds e.g.~to \cite[Eq.~(2.13)]{CI}. As an immediate consequence of the above proposition we get the following (the proof is omitted since it is trivial):

\begin{Corollary}  Denoting by $\nabla \theta (\xi,z)$ the vector $\big(\nabla \theta_1 (\xi,z),\cdots,  \nabla \theta_d (\xi,z) \big)$ in $\bbC^d$, for  any $a\in\bbR^d$ we have
\begin{equation}           
    \sigma(\o)a
        =2 \sum_{z\in\mathcal N_*}\bbE\Big[c_{0,z}\big( z\cdot a+\nabla\theta(\cdot,z)\cdot a\big)\Big]z
\end{equation}
so that in particular
\begin{equation}\label{carlton}        
    a\cdot\sigma(\o)a
        =2\sum_{z\in\mathcal N_*}\bbE\Big[c_{0,z}\big( z\cdot a+\nabla\theta(\cdot,z)\cdot a\big)\Big](z\cdot a)\,.
\end{equation}
In the nearest-neighbour case these reduce to
\begin{equation}           
    \sigma(\o)a
        =2\sum_{j=1}^d\bbE\Big[c_{0,e_j}\big( a_j+\nabla\theta(\cdot,e_j)\cdot a\big)\Big]e_j
\end{equation}
and
\begin{equation}\label{rododendro}           
    a\cdot\sigma(\o)a
        =2\sum_{j=1}^d\bbE\Big[c_{0,e_j}\big( a_j+\nabla\theta(\cdot,e_j)\cdot a\big)\Big]a_j\,.
\end{equation}
\end{Corollary}
We point out that formula \eqref{rododendro} corresponds e.g.~to the first identity in \cite[(3.6)]{CI}. On the other hand, the second identity in \cite[(3.6)]{CI}
does not hold  in the framework of the present paper. One needs an additional term:

\begin{Proposition}
    For all $a\in\bbR^d$ we have
    \begin{equation}\label{pavone1}
        a\cdot\sigma(\o)a=-2i\o\bbE[|\theta\cdot a|^2]+2\sum_{z\in\mathcal N_*}\bbE\big[c_{0,z}|\nabla\theta(\cdot,z)\cdot a+z\cdot a|^2\big]\,.
    \end{equation}
    In the nearest-neighbour case this formula reduces to
    \begin{equation}\label{pavone2}
        a\cdot\sigma(\o)a
            =-2i\o\bbE[|\theta\cdot a|^2]+2\sum_{j=1}^d\bbE\big[c_{0,e_j}|\nabla\theta(\cdot,e_j)\cdot a+ a_j|^2\big]\,.
    \end{equation}    
\end{Proposition}
We recall that $\theta(\xi)$ and $\nabla \theta(\xi,z)$ are vectors in $\bbC^d$. We stress that \eqref{pavone1} and \eqref{pavone2} split the real and imaginary parts of  $a\cdot\sigma(\o)a$.
\begin{proof}
    By \eqref{carlton} one can write
    \begin{align}\label{blando}
        a\cdot\sigma(\o)a
            &=2\sum_{z\in\mathcal N_*}\bbE\big[c_{0,z} (z\cdot a)^2\big]+2\sum_{z\in\mathcal N_*}\bbE\big[c_{0,z}(\nabla \theta(\cdot,z)\cdot a)(z\cdot a)\big]\,.
    \end{align}
    We can rewrite the second summand as
    \begin{align}\label{lemming}
    2\sum_{z\in\mathcal N_*}&\bbE\big[c_{0,z}(\nabla \theta(\cdot,z)\cdot a)(z\cdot a)\big]\nonumber\\
        &=4\sum_{z\in\mathcal N_*}\bbE\big[c_{0,z}(\nabla \theta(\cdot,z)\cdot a)(z\cdot a)\big]+2\bbE\big[(\theta\cdot a)(\gamma\cdot a)\big]\,,
    \end{align}
    where we have used the fact that $\bbE[c_{0,z}(\theta(\tau_z\cdot)\cdot a)(z\cdot a)]=\bbE[c_{0,-z}(\theta\cdot a)(z\cdot a)]$ and the definition $\gamma(\xi)=\sum_{z\in\mathbb N_*}(c_{0,z}(\xi)-c_{0,-z}(\xi))z$.
    Let us analyze the expectation $\bbE[(\theta\cdot a)(\gamma\cdot a)]$.
    Since $\gamma=(i\o-\bbL)\theta$ and $\g$ is real (thus implying that $\gamma=(-i\o-\bbL)\bar\theta$)  we have
    \begin{equation}\label{ciriola}
        \bbE\big[(\theta\cdot a)(\gamma\cdot a)\big]
            =- i\o\bbE\big[ |\theta\cdot a|^2\big]+\bbE\big[(\theta\cdot a)(-\bbL (\bar\theta\cdot a))\big]\,.
    \end{equation}
    We write $\bbL_R$ for the original operator  $\bbL$ in Proposition~\ref{zap} to distinguish it from its complex extension, denoted by $\bbL$ above.
    As $\theta_k \in \cD(-\bbL)$, we have 
    that $\Re(\theta_k)$ and $\Im(\theta_k)$ belong to $\cD(-\bbL_R)$ inside $ L^2_R(\bbP)$. On the other hand,  by spectral calculus, $\cD(-\bbL_R)\subset \cD(\sqrt{-\bbL_R})$. By the parallelogram rule, from \eqref{def_dir_form}  and \eqref{istria} in Proposition~\ref{zap} we derive
\be \label{castagna}
\la \sqrt{-\bbL_R}f , \sqrt{-\bbL_R}g\ra_{L^2(\bbP)}= 
	\frac{1}{2} \sum _{z\in\cN} 
	\bbE\big[ c_{0,z}(\xi) \left( f (\t_z \xi) - f(\xi)\right)\, \left( g (\t_z \xi) - g(\xi)\right) \big] 
\en
for any $f,g\in \cD(\sqrt{-\bbL_R})\subset L^2_R(\bbP)$. Applying the above formula to $f,g$ of the form $\Re(\theta\cdot a)$ and $\Im(\theta\cdot a)$ we finally get 
        \begin{align}\label{teddy}
        \bbE[(\theta\cdot a)(-\bbL (\bar\theta\cdot a))]
            =\frac12\sum_{z\in\mathcal N}\bbE\Big[c_{0,z}  \big|\nabla \theta(\cdot,z)\cdot a\big|^2 \Big]
            =\sum_{z\in\mathcal N_*}\bbE\Big[c_{0,z}\big|\nabla \theta(\cdot,z)\cdot a\big|^2\Big]\,.
    \end{align}
    Putting together \eqref{lemming}, \eqref{ciriola} and \eqref{teddy}  back into \eqref{blando}, we have arrived to
    \begin{align}
    a\cdot\sigma(\o)a
        &= 2\sum_{z\in\mathcal N_*}\bbE\Big[c_{0,z} (z\cdot a)^2\Big]
        +4\sum_{z\in\mathcal N_*}\bbE\Big[c_{0,z}(\nabla \theta(\cdot,z)\cdot a)(z\cdot a)\Big] \nonumber\\
        &-2i\o\bbE[|\theta\cdot a|^2]
        +2\sum_{z\in\mathcal N_*}\bbE\Big[c_{0,z}\big|\nabla \theta(\cdot,z)\cdot a\big|^2\Big]\nonumber
    \end{align}    
    which is equivalent to \eqref{pavone1}. Equation \eqref{pavone2} follows from \eqref{pavone1}.
 \end{proof}

\section{Preliminary results}\label{sec_hom}
In this section we introduce some concepts (including the 2-scale convergence) and technical tools to get the relevant homogenization results of  the next section. Our approach  is based on  \cite{F_hom} and \cite{F_resistor}, where the corresponding proofs can be found. The simpler context here allows  for a significant simplification of those derivations. We point out that \cite{F_hom} treats much more general models and our  conductance model  can be retrieved 
by setting there $\bbG:=\bbZ^d$, $r_{x,y}(\o):=c_{x,y}(\o)$, $n_x(\o)=1$, $\hat\o:=\bbZ^d$, $\o$ denoting the environment in \cite{F_hom}. The moment conditions in \cite{F_hom} and \cite{F_resistor} then reduce to $\bbE[c_{x,y}]<+\infty$. In this section we will state explicitly where we use the stronger assumption $c_{0,z}\in L^2(\bbP)$ for all $z\in \cN$, for the rest we will just assume without further mention that $c_{0,z}\in L^1(\bbP)$. 

Although  the results in \cite{F_hom} refer to  real functions, their extension  to complex functions can be obtained immediately  by distinguishing the contribution of the  real and imaginary parts.

We define $\nu$ as the Radon measure on $\O \times \cN$ such that 
\begin{equation}\label{labirinto}
\int_{\O\times \cN} {\rm d} \nu (\xi, z) g (\xi, z) 
	= \int_\O  {\rm d} \bbP (\xi) \sum_{z\in \cN}  c_{0,z}(\xi) g( \xi, z) 
\end{equation} 
 for any nonnegative   measurable function $g:\Omega\times\cN\to\bbR$. Since $ c_{0,z}\in L^1(\bbP)$ for all $z\in \cN$  by hypothesis, $\nu$ has finite total mass.

 Recall that, given a  function $u:\O\to \bbC$,   the gradient $\nabla u: \O \times \cN \to \bbC$  of $u$ has been defined in \eqref{cantone} as  $ \nabla u (\xi, z):= u (\t_z \xi)-u (\xi)$.
 If $u $ is bounded and measurable, then $\nabla u \in L^2(\nu)$.
The subspace of \emph{potential forms} $L^2_{\rm pot} (\nu)$ is defined as   the  closure in $L^2(\nu)$ of gradient functions:
 \[ L^2_{\rm pot} (\nu) :=\overline{ \{ \nabla u\,:\, u \text{ is  bounded and measurable} \}}\,.
 \]

\begin{Definition}\label{def_div}
Given  $v:\O\times \cN \to \bbC$ its divergence  $\dive v: \O\to\bbC$ is defined as 
\begin{equation}\label{emma_0}
\dive v(\xi):= \sum_{z\in \cN}c_{0,z}(\xi) ( v(\xi,z)-  v(\t_z \xi, -z) )\,.
\end{equation}
\end{Definition}
Trivially, if  $v=v'$  $\nu$--a.s., then  $\dive v = \dive v'$ $\bbP$--a.s.. In particular, the definition of $\dive v$ extends to $L^p$--spaces. By Cauchy--Schwarz inequality,  for any  $v \in L^2(\nu)$ it holds $\dive v\in L^1(\bbP)$. 
By using the stationarity of $\bbP$ one easily gets the following useful integration by parts rule:
\begin{Lemma}\label{faticaccia} Consider measurable functions  $u: \O\to \bbC $ and $v:\O\times \cN \to \bbC$.  Suppose that at least one of the following conditions is satisfied:
\begin{itemize}
\item[(i)]   $v \in L^2(\nu)$ and $u$ is bounded;
\item[(ii)] $v$ is bounded, $u\in L^2(\bbP)$ and $c_{0,z}\in L^2(\bbP)$ for all $z\in \cN$.
\end{itemize}
Then 
\begin{equation}\label{italia}
\int _\O {\rm d} \bbP(\xi)   \dive v(\xi)  u (\xi)= - \int_{\O\times \cN} {\rm d} \nu(\xi, z) v( \xi, z) \nabla u (\xi, z) \,.
\end{equation}
The above integrands are integrable and in particular the integrals are well defined.
\end{Lemma}
\begin{proof}
The case of item (i) corresponds to   \cite[Lemma~8.3]{F_hom} (the verification here is trivial). For item (ii) the derivation is again trivial since our assumptions imply that $c_{0,z}(\xi) v(\xi,z)$, $c_{0,z}(\xi) v(\t_z \xi, -z) $, $c_{0,z}(\xi) u(\xi)$ and $c_{0,z}(\xi) u (\t_z \xi) $ belong to $L^1(\bbP)$. Moreover, by the stationarity of $\bbP$ and since $c_{0,z}(\xi)=c_{-z,0} (\t_z \xi)$ by the covariant relation \eqref{covariante}, we have
\[
\int_\O\rmd \bbP(\xi) c_{0,z}(\xi) v(\t_z \xi, -z)u(\xi)= \int_\O \rmd \bbP(\xi) c_{0,-z}(\xi) v(\xi, -z) u (\t_{-z}\xi)\,.
\]
This identity allows to check \eqref{italia} as $\cN$ is symmetric.
\end{proof}
%
%

\subsection{Set $\O_{\rm typ}$ of typical environments}\label{sec_typ_env}  
In the rest of the paper, we will work on a full-measure restriction $\O_{\rm typ}$ of the set $\Omega$, that we will call the set of typical environments. Our definition of $\O_{\rm typ}$ is related to the one appearing in \cite[Section 12]{F_hom}, since we want to recover some properties for $\o\in\O_{\rm typ}$ that have been proved in \cite{F_hom}. On the other hand, the context of the present paper allows for a simplification of the definition given in \cite{F_hom}.  

\smallskip

$\O_{\rm typ}$ is defined as the intersection of several sets (see Definition \ref{def_om_typ} below): for their definition we need to introduce a few more objects.
Recall that, as proved in \cite{F_hom} using the assumption that $L^2(\bbP)$ is separable, also   $L^2(\nu)$ is separable. 


\smallskip

\noindent
{\bf The functional sets $\cG_2,\cH_2,\cH_3 $}.  We fix a countable set $\cG_2$  of  bounded real measurable functions 
$g: \O\to \bbR$ such that  the set $\{ \nabla g\,:\, g \in \cG_2\}$, thought of  in $L^2_R(\nu)$,  is  dense in $L^2_{\rm pot}(\nu) \cap L^2_R(\nu)$.  At the cost of adding a countable family of bounded functions (recall that $L^2_R(\bbP)$ is separable by assumption), $\cG_2$ is dense in $L^2_R(\bbP)$. We  define $\cH_2$ as the set of measurable functions  $h: \O \times \cN \to \bbR$ such that $h=\nabla g$ for some $g\in \cG_2$. We define $\cH_3$ as  the set of measurable functions $h: \O \times \cN \to \bbR$ such that $h(\xi,z)= g(\t_z \xi) z_i$ for some $g\in \cG_2$ and some direction $i=1,\dots, d$. Note that  $\|h\|_{L^2(\nu)} <+\infty$ for all $h \in \cH_2\cup \cH_3$  
since functions $g\in \cG_2$ are bounded.

The functional sets $\cG_2,\cH_2,\cH_3 $ appear also in \cite[Section 12]{F_hom} and we have kept the same name to allow a comparison in some arguments appearing in the rest (although this implies that  the numbering starts from  $2$).

\smallskip

\noindent
{\bf The functional set $\cH_*$}. 
We fix  a countable set  $\cH_*$ of bounded measurable real functions dense in  $L^2_R(\nu) $ (the  existence follows from the separability of $L^2_R(\nu)$).




\begin{Definition}\label{ometto}
Given  a   function  $b :\O \times \cN \to \bbR$ we define $\tilde b: \O \times \cN \to \bbR$ as 
$
\tilde b (\xi, z):=
 b (\t_{z} \xi, -z)$. Moreover, we define $\hat b : \O \to \bbR$ as $\hat b (\xi) :=\sum_{z\in \cN} c_{0,z}(\xi) b(\xi, z) $.
\end{Definition}

For the next defintion we need to assume that $c_{0,z}\in L^2(\bbP)$.
\begin{Definition}[Functional  set $\cG$]  \label{gonars}
We take as $\cG$ the  countable set of measurable functions on $\O$ given by the union of the following subsets of $L^2_R(\bbP)$:  $\{1\}$, $\cG_2$, $c_{0,z}$ for $z\in \cN$,  $\{ \dive h\,:\,h\in \cH_*\}$,   and   $\{h_i \,:\, i=1,..,d \text{ with } b \in\cH_*\}$, where  
$h_i(\xi):=\sum_{z\in \cN} c_{0,z}(\xi)  z_i \tilde b(\xi,z)$. 
\end{Definition}
Note that $\cG$ is dense in $L^2_R(\bbP)$ since it contains $\cG_2$ which is dense. 
Note also  that, since $c_{0,z}\in L^2_R(\bbP)$ and all $h\in\cH_*$ are bounded, we have $\dive h \in L^2_R(\bbP)$ for any $h\in \cH_*$ and  similarly $h_i\in L^2_R(\bbP)$ with $h_i$ as in Definition~\ref{gonars}.

\begin{Definition}[Functional  set $\cH$] We  take as  $\cH$  the    countable set  of measurable functions on $\O\times \cN$   given by the union of  the following subsets of $L^2_R(\nu)$:  $\cH_2$, $\cH_3$, $\cH_*$, 
$\{\tilde b : b \in  \cH_2\cup \cH_3\cup \cH_*\}$.
\end{Definition}
Note that $\cH$ is dense  $L^2_R(\nu)$ since it contains $\cH_*$ and that $\tilde b \in \cH $ for all $b\in \cH$ (since $\tilde{\tilde b}= b$).



To define the set of typical environments we need a useful consequence of  the ergodic theorem stated 
 e.g.~in \cite{F_hom} (see \cite{F_x_Francis} for a proof and generalizations): 
\begin{Proposition}\label{prop_ergodico} \emph{(cf.~\cite[Proposition~3.1]{F_hom})}
 Let  $g: \O \to \bbR$ be a measurable function with $\|g\|_{L^1(\bbP)}<+\infty$. Then there exists a translation invariant   measurable subset $\cA[g]\subset \O$  such that $\bbP(\cA[g])=1$ and such that,  for any $\xi\in \cA[g]$ and any  $\varphi \in C_c (\bbR^d)$, it holds
\begin{equation}\label{eq_limitone_base}
\lim_{N\to +\infty}\frac{1}{N^d} \sum _{x\in N^{-1}\bbZ^d} \varphi(x) g(\t_{Nx} \xi) =
\int _{\bbR^d} {\rm d}x\,\varphi (x) \, \bbE[g]\,. 
\end{equation}
\end{Proposition}

\begin{Definition}[Set of typical environments $\O_{\rm typ}$]  \label{def_om_typ}
We define the set $\O_{\rm typ}$ as the intersection of the following sets: $\cA[g g' ]$ for all $g,g'\in \cG$; 
 $\cA[\widehat{b b'}]$ for all $b,b'\in \cH$.
\end{Definition}
Note that the above  definition is well posed. Indeed,  for $g,g' \in \cG$ we have  $gg'\in L^1(\bbP)$ since $\cG\subset L^2 (\bbP)$; while for $b,b'\in \cH$ we have  $\widehat{b b'}\in L^1(\bbP)$ as can be easily derived from Cauchy--Schwarz inequality,  the identity $\widehat{b b'}(\xi) =\sum _{z\in \cN} c_{0,z}(\xi) b(\xi,z) b'(\xi,z)$ and the property $\cH\subset L^2(\nu)$.

Note that $\O_{\rm typ}$ is  measurable, translation invariant and has $\bbP$--probability equal to $1$ being a countable intersection of sets  fulfilling the same properties.
For later use, we point out that the inclusions 
$\O_{\rm typ}\subset \cA[c_{0,z}]$ and $\O_{\rm typ}\subset \cA[c_{0,z}^2]$ in Definition~\ref{def_om_typ}, imply respectively that  for all $\xi\in \O_{\rm typ}$   it holds
\begin{align}
& \lim_{N\to +\infty}\frac{1}{N^d} \sum _{x\in [0,N)^d\cap \bbZ^d}  c_{x,x+z}(\xi)  =\bbE[ c_{0,z} ] \qquad \forall z\in \cN\,,\label{req1} \\
& \lim_{N\to +\infty}\frac{1}{N^d} \sum _{x\in [0,N)^d\cap \bbZ^d}  c^2_{x,x+z}(\xi)  =\bbE[ c^2_{0,z} ] \qquad \forall z\in \cN\,.\label{req2} 
\end{align}

\begin{Remark}The functions in $\cG$ and $\cH$ will play the role of test functions. We just need to test against real functions.  Similarly in what follows $C^\infty_c(\L)$ and $C_c(\L)$ will consist of real functions. 
\end{Remark}

\subsection{Measures $\mu_\e$ and $\nu_\e^\xi$}\label{baglioni}
Let $\L:=(0,1)^d\subset \bbR^d$. Set 
\be\label{mango}   
\e:= \frac{1}{N}\,,\qquad \L_\e :=\L \cap ( \e \bbZ^d )\,.
\en
We define 
$\mu_\e$, $\nu_\e^\xi$  as the atomic measures on $\L$ and $\L\times \cN$, respectively, given by 
\be\label{onde_alpha}
\mu_\e :=\e^d \sum_{x\in \L_\e} \d_x  \,, \qquad \;\;
\nu_\e^\xi:=\e^d \sum_{x\in \L_\e} \sum_{\substack{z\in \cN:\\ x+\e z\in \L_\e }} c_{x/\e,x/\e+ z} (\xi)\d_{(x,z)}\,.
\en
We trivially point out that $\mu_\e$ does not depend on $\xi$.


An immediate consequence of Proposition~\ref{prop_ergodico} is the following fact:
\begin{Corollary}\label{cor_ergodico} 
 Let  $g: \O \to \bbR$ be a measurable function with $\|g\|_{L^1(\bbP)}<+\infty$. Then for any $\xi\in \cA[g]$ and any  $\varphi \in C_c (\L)$ it holds
\begin{equation}\label{eq_limitone}
\lim_{\e\da 0} \int _\L {\rm d}  \mu_\e  (x)  \varphi (x )g(\t_{x/\e} \xi)=
\int _{\L} {\rm d}x\,\varphi (x) \, \bbE[g]\,. 
\end{equation}
\end{Corollary}

\subsection{2-scale convergence  in $L^2( \mu_\e)$}\label{sec_2scale1}
We recall the definition of 2-scale convergence   in our context for  $L^2( \mu_\e)$:
\begin{Definition}\label{priscilla}
Fix   $\z\in \O_{\rm typ}$, an $\e$--parametrized  family of functions  $v_\e \in L^2(\L,  \mu_\e)$ and  a function $v \in L^2 \bigl(\L\times \O, {\rm d}x \times \bbP \bigr)$.
\\
$\bullet$ We say that \emph{$v_\e$ is weakly 2-scale convergent to $v$}, and write 
$v_\e \stackrel{2}{\rightharpoonup} v$, 
if the family $\{v_\e\}$ is asymptotically bounded, i.e.
$
 \limsup_{\e\downarrow 0}  \|v_\e\|_{L^2(\L,\mu_\e)}<+\infty$, 
 and 
\begin{equation}\label{rabarbaro}
\lim _{\e\downarrow 0} \int_\L  {\rm d} \mu _\e (x)  v_\e (x) \varphi (x) g ( \t _{x/\e}\zeta ) 
	=\int_\O {\rm d}\bbP(\xi)\int_\L {\rm d}x\,  v(x, \xi) \varphi (x) g (\xi)  \,,
\end{equation}
for any $\varphi \in  C_c (\L)$ and any  $g \in  \cG$. \\
$\bullet$
 We say that \emph{$v_\e$ is strongly 2-scale convergent to $v$}, and write 
$v_\e \stackrel{2}{\to} v$, 
if the family $\{v_\e\}$ is asymptotically bounded and 
\begin{equation}\label{gingsen}
\lim_{\e\downarrow 0} \int_\L {\rm d}\mu_\e(x)   v_\e(x) u_\e(x) = \int _\O d\bbP (\xi)\int_\L  {\rm d}x\,   v(x, \xi) u(x,\xi)  
\end{equation}
whenever $u_\e \stackrel{2}{\rightharpoonup} u$. 
\end{Definition}
Although we just write $v_\e \stackrel{2}{\rightharpoonup} v$ and $v_\e \stackrel{2}{\to} v$, the above convergences depend on the environment $\z$ which indeed appears in the l.h.s. of \eqref{rabarbaro}.

\begin{Lemma}\label{compatto1} \emph{(cf.~\cite[Lemma~10.3]{F_resistor})}
 Let $\z\in \O_{\rm typ}$.    Then, given an asymptotically bounded family of functions $v_\e\in L^2 (\L, \mu_\e )$,  there exists a vanishing sequence $(\e_k)$   such that
 $ v_{\e_k}
  \stackrel{2}{\rightharpoonup}   v $ for some 
 $  v \in L^2(\L\times \O,  {\rm d}x \times \bbP )$ with $\|  v\|_{   L^2(\L \times \O,{\rm d}x \times \bbP )}\leq \limsup_{\e\da 0} \|v_\e\|_{L^2( \mu_\e)}$.
 \end{Lemma}
 
We point out that the sequence $(\e_k)$ depends on $\z$.
The proof of the above lemma is an adaptation of the proof given in \cite[Appendix~F.1]{F_hom} to the present setting. 
  It relies only on the property that $\cG$ is countable and dense in $L^2_R(\bbP)$ and that $\O_{\rm typ}\subset \cA[g g' ]$ for all $g,g'\in \cG$.

\subsection{2-scale convergence in $L^2( \nu _\e^{\z})$}\label{sec_2scale2}
Recall the definition of the measure $\nu$ given in \eqref{labirinto}.
 We recall the definition of 2-scale convergence   in our context for  $L^2( \nu^\z_\e)$:
 \begin{Definition}\label{chioccia}
Given  $\z\in \O_{\rm typ}$, an $\e$--parametrized  family of functions $w_\e \in L^2(\L\times \cN, \nu _\e^\z)$ and  a function  $w \in L^2 \bigl(\L\times \O\times \cN\,, \rmd x \times \rmd \nu\bigr)$, we say that \emph{$w_\e$ is weakly 2-scale convergent to $w$}, and write $  w_\e \stackrel{2}{\rightharpoonup}w   $,  if   $\{w_\e\}$ is asymptotically bounded
 in $L^2(\L\times \cN, \nu _\e^\z)$, i.e. $
 \limsup_{\e \da 0} \|  w_\e\|_{L^2(\L\times \cN, \nu _\e^\z)}<+\infty$,
   and
\be\label{yelena}
 \lim _{\e\downarrow 0} \int _{\L\times \cN}  \rmd \nu_\e^\z (x,z) w_\e (x,z ) \varphi (x) b ( \t _{x/\e} \z ,z )\\
=\int_{\L}  \rmd x \int _{\O\times \cN} \rmd \nu (\xi, z) w(x, \xi,z ) \varphi (x) b (\xi,z )  \,,
\en
for any    $\varphi \in C_c(\L)$  and  any  $b \in \cH $. \end{Definition}
Also the above 2-scale convergence $  w_\e \stackrel{2}{\rightharpoonup}w   $ depends on the environment $\z$, although the latter does not appear in the notation.

\begin{Lemma}\label{compatto2} \emph{(cf.~\cite[Lemma~10.5]{F_resistor})}  Let $\z \in \O_{\rm typ}$.  Then, given an asymptotically bounded family of functions $w_\e\in L^2 (\L, \nu_\e^\z )$,  there exists a vanishing sequence $(\e_k)$   such that  $w_{\e_k} \stackrel{2}{\rightharpoonup} w$ for some 
 $ w \in L^2(  \L \times \O\times\cN\,,\, dx \times \nu )$ with  $\|  w\|_{   L^2(\L\times \O\times \cN\,,\, \rmd x \times \nu )}\leq \limsup_{\e\da 0} \|w_\e\|_{L^2(\L\times \cN, \nu_\e^\z)}$.
\end{Lemma}
We point out that in the  lemma above the sequence $(\e_k)$ depends on $\z$. 
The proof of the above lemma is an adaptation of the proof given in \cite[Appendix~F.2]{F_hom} to the present setting. 
  It relies only on the assumption  that $\cH$ is countable and dense in $L^2_R(\nu)$ and that $\O_{\rm typ}\subset \cA[\widehat{b b'}]$ for all $b,b'\in \cH$.

\subsection{$\e$--gradient}\label{eps_gradient}
While the gradient presented in \eqref{cantone} is associated to functions acting on elements of $\Omega$, we can also introduce  a microscopic spatial gradient (that we call \emph{$\e$--gradient}) associated to functions $v:\,\bbR^d\to\bbR$ as
\begin{align}\label{gradientello}
\nabla_\e v(x,z)
	:= \frac{v(x+ \e z)- v(x)}{\e} \qquad\mbox{ for $x\in \e \bbZ^d, z\in \cN$}\,.
\end{align}
One can easily check the following Leibniz rule:
\begin{equation}\label{leibniz} \nabla _\e  (fg)(x,z)
   =\nabla _\e  f (x, z ) g (x )+ f (x+\e z ) \nabla _\e g  ( x, z)\,.
\end{equation}
We further introduce the equivalent of the measures in \eqref{onde_alpha}  defined on the whole space:
\be\label{onde_beta}
\bar \mu_\e :=\e^d \sum_{x\in \e \bbZ^d } \d_x  \,, \qquad \;\;
\bar \nu_\e^\xi:=\e^d \sum_{x\in \e\bbZ^d } \sum_{\substack{z\in \cN }} c_{x/\e,x/\e+ z}  (\xi)\d_{(x,z)}\,.
\en
 We have the following results:
 \begin{Lemma}\emph{(cf. \cite[Lemma~11.3]{F_hom})}  \label{gattonaZ}   
 Let $b: \O \times \bbZ^d \to\bbC$ and  let  $\varphi, \psi:\bbR^d \to \bbC $  be
   functions with bounded support. Then, for each $\xi\in\Omega$ and $\e>0$,  it holds
  \begin{align}
     \int \rmd \bar\nu_\e^\xi  (x, z) \varphi ( x) \psi (  x+\e z ) b(\t_{x/\e}\xi, z)
     	&= \int \rmd \bar \nu^\xi_\e (x, z) \psi (x) \varphi (  x+\e z  ) \tilde{b}(\t_{x/\e}\xi, z)\,, \label{micio2}
\\
    \int \rmd \bar \nu_\e^{ \xi} (x, z) \nabla_\e  \varphi   ( x,z) \psi (   x+\e z ) b(\t_{x/\e}\xi, z) 
    &=  - \int \rmd\bar  \nu_\e^{ \xi}  (x, z) 
  \nabla _\e \varphi  (x,z) 
  \psi (x) \tilde{b}(\t_{x/\e}\xi, z)\,.\label{micio3}
   \end{align}
\end{Lemma}
\begin{Lemma}\emph{(cf.~\cite[Lemma~11.7]{F_hom})}  \label{cane_ciak}
Let  $b: \O\times \bbZ^d \to \bbC$ and let $u:\bbR^d\to \bbC$   be a function with bounded support.   Then, for each $\xi\in\Omega$ and  $\e>0$,  it holds
\begin{equation}\label{sea}
\int \rmd \bar \mu_\e (x) u(x) \dive b (\t_{x/\e} \xi) = - \e \int \rmd \bar\nu _\e^\xi (x,z) \nabla_\e u(x,z) b ( \t_{x/\e} \xi, z) \,.
\end{equation}
\end{Lemma}
Note that the integrals in \eqref{micio2}, \eqref{micio3} and \eqref{sea} are actually sums with a finite number of addenda.
Our context is much simpler than the general one in \cite{F_hom}  and a direct check of Lemma~\ref{gattonaZ} and Lemma~\ref{cane_ciak} (based on integration by parts) is very simple.

\begin{Remark}\label{passaggio}
Since $\cN$ is a bounded set, if $\varphi$ or  $\psi$ 
have support in $\{x\in \L\,:\, d(x, \partial \L)\geq\d\}$ for some $\d>0$, then \eqref{micio2} and   \eqref{micio3}  hold with
 $\nu   _\e ^\xi$ instead of $\bar \nu   _\e ^\xi$ for all $\e \leq \e_*( \cN, \d)$ for some $\e_*(\cN,\d)>0$. Similarly, if $u$ has support 
 $\{x\in \L\,:\, d(x, \partial \L)\geq \d\}$ for some $\d>0$, then \eqref{sea} holds with 
 $\mu_\e$ instead of $\bar \mu_\e^\xi$  and $\nu   _\e ^\xi$ instead of $\bar \nu   _\e ^\xi$ for all $\e \leq \e_*( \cN, \d)$.
\end{Remark}
For later use we recall the following result: 
\begin{Lemma}\emph{(cf.~ \cite[Lemma 19.2]{F_hom})} \label{pesciolino} Given  $\zeta\in \O_{\rm typ}$ and $\varphi\in C_c^2 (\L)$ it holds
\be\label{agosto}
\lim_{\e \da 0} \int \rmd\nu_\e^\zeta (u,z) [\nabla_\e \varphi (u,z)-\nabla \varphi (u) \cdot z]^2=0\,.
\en 
\end{Lemma}
In \eqref{agosto},  $\nabla \varphi(u)$ is the standard gradient. In the present context, the derivation of this lemma is trivial and we give it for completeness.
\begin{proof}By Taylor expansion, given $u\in \L_\e$ and $z\in \cN$, we have $\big|\varphi(u+ \e z) -\varphi (u)- \e \nabla \varphi (u) \cdot z\big| \leq c(\cN,\varphi) \e^2$ for some positive constant $c(\cN,\varphi)$ determined only by $\cN,\varphi$.
Therefore, \eqref{agosto} holds whenever $\limsup_{\e \da 0} \int \rmd\nu_\e^\zeta (u,z)<+\infty$. We have for $\e $ small that 
\begin{equation*}
\begin{split}  \int_{\L\times \cN} \rmd\nu_\e^\zeta (u,z)&= \e^d \sum_{x\in \L_\e} \sum_{\substack{z\in \cN:\\ x+\e z\in \L_\e }} c_{x/\e,x/\e+ z} (\z)\leq  \e^d \sum_{x\in \L_\e} \sum_{z\in \cN}c_{x/\e,x/\e+ z} (\z)
\,.
\end{split}
\end{equation*}
Since $\bbE[c_{0,z}]<+\infty$ for all $z\in \cN$, the claim then follows from property \eqref{req1}, which is  fulfilled  by $\z\in \O_{\rm typ}$.
 \end{proof}

\section{Proof of Theorem \ref{teo1bis} (I): 2-scale convergence}\label{sec_dim_I}
In this section we come back to the torus $\bbT^d_N$ and set $\e:=1/N$. Without loss we restrict to $N>2 \|\cN\|_\infty$. We make a simple but very relevant observation. Given $x\in \L_\e$  and $z\in \cN$ with $x+\e z\in \L_\e$, it holds
\be
c_{x/\e, x/\e+z}(\xi) = c^{(N)}_{x/\e, x/\e+z}(\xi)\,.
\en
Roughly put, given two points in the bulk we do not see the effect of the environment periodization. This observation will be frequently used in the second part of this section and allows to use the concepts of the previous section.

We recall that $\mfm_N$ is the uniform probability  measure on $\bbT^d_N$. We also recall that, since we deal with complex functions,  we think of  $L^2(\mfm_N)$ as the Hilbert space of complex functions on $\bbT^d_N$ with scalar product
\be
\la f, g\ra _{L^2(\mfm_N)}:= \int _{\bbT^d_N} \rmd \mfm_N(x) \bar f (x) g(x)= N^{-d} \sum_{x\in \bbT^d_N} \bar f(x) g(x)\,.
\en
Recall  (cf.~\eqref{cancellami})  that we associated to $\xi\in\Omega$ and $N$ a  symmetric operator on $ L^2(\bbT^d_N,\mathfrak{m}_N)$  given by 
$
\cL^\xi_N f(x)=\sum_{z\in \cN}c^{(N)}_{ x,  x+z}(\xi) (f(x+z)-f(x))$ for $ x\in\bbT^d_N$.
It is simple to check that
\be\label{sinto} 
\begin{split}
\la -\cL_N^\xi f,  g\ra _{L^2(\mathfrak{m}_N)}&=
\la f, -\cL_N^\xi g\ra _{L^2(\mathfrak{m}_N)}
	\\
    &=\frac{1}{2} N^{-d}\sum _{x\in \bbT^d_N} \sum_{y: \,y\sim x} c ^{(N)}_{x,y}(\xi)  \left( \bar f(y)-\bar f(x) \right) \left( g(y)-g(x) \right)\,,
\end{split}
\en
where now $y\sim x$ means that $y=x+z$ with $z\in \cN$. 

Let $\nu_N^\xi$ be the atomic measure on $\bbT^d_N\times \cN$ given by 
\be
\nu_N^\xi:=N^{-d} \sum_{x\in \bbT^d_N} \sum_{z\in \cN} c^{(N)}_{x,x+z}(\xi)\d_{(x,z)}\,.
\en
The space $L^2(\nu_N^\xi)$ is given by  complex functions and the associated scalar product is 
\be
\la F, G\ra_{L^2( \nu_N^\xi)}:=N^{-d} \sum _{(x,z)\in  \bbT^d_N\times \cN} c^{(N)}_{x,x+z}(\xi) \bar F (x,z)  G(x,z) \qquad
F,G: \bbT^d_N\times \cN\to \bbC\,. 
\en

Given $f:\bbT^d_N\to \bbC$ and $z\in\mathcal N$,  we set  
\be\label{grillo} \nabla f(x,z) := f(x+z)-f(x) \,.
\en
Then  \eqref{sinto} is equivalent to
\be\label{rom}
\la -\cL_N^\xi f,  g\ra _{L^2(\mathfrak{m}_N)}=
\la f, -\cL_N^\xi g\ra _{L^2(\mathfrak{m}_N)}=\frac{1}{2}\la \nabla f, \nabla g \ra _{L^2(\nu^\xi_N)}\,.
\en

\begin{Remark} Overall we deal with three types of gradient: the $\e$--gradient $\nabla _\e f (x, z)=\e^{-1}\big( f(x+\e z) -f (x)\big)$ with $x\in \e\bbZ^d$ introduced in \eqref{gradientello}, the toroidal spatial gradient $\nabla f(x,z):= f(x+z)-f(x)$  with $x\in \bbT^d_N$ introduced in \eqref{grillo} and the environment  gradient  $\nabla f(\xi,z):= f(\t_z \xi) -f (\xi)$ with $\xi\in \O$ introduced in \eqref{cantone}.
\end{Remark}
In what follows, given a complex function $f$, we write $f_R$ and $f_I$ for its real and imaginary part, respectively, so that $f=f_R+if_I$. One can easily check that, for any $f:\bbT^d_N\to \bbC$, it holds
\be\label{romix} \frac{1}{2}\la \nabla f, \nabla f \ra _{L^2(\nu^\xi_N)}=
\frac{1}{2}\la \nabla f_R, \nabla f_R \ra _{L^2(\nu^\xi_N)}+ \frac{1}{2}\la \nabla f_I, \nabla f_I \ra _{L^2(\nu^\xi_N)}\,.
\en


\smallskip

 Recall the definition \eqref{def_gamma_N}  of the local drift  $\g^\xi_N$. We fix a vector $v\not=0$ and denote the projection of the local drift in direction $v$ with $ \gamma_{N,v}^\xi:\, \bbT^d_N \to\bbR$, i.e.
 \be 
  \gamma_{N,v}^\xi(x)
 	:= \g_N^\xi(x) \cdot v
 	=\sum_{z\in \cN} c^{(N)}_{x,x+z}(\xi) (z\cdot v)\,. 
 \en
 Recall  definition~\eqref{def_gammaV} of $\g(\xi)$. Analogously to what done for $\g^\xi_{N,v}(x)$, we define the function
 $\gamma_v:\,\Omega\to\bbR$ as 
 \be \gamma_v(\xi):=\gamma(\xi)\cdot v =\sum_{z\in \cN} c_{0,z}(\xi) (z\cdot v)\,.\en
 For   later use we stress that, trivially, $ \gamma_{N,v}^\xi$ and $\gamma_v$ are real functions.
 
We introduce the  function $\theta^\xi_N:\,\bbT^d_N\to\rosso{\bbR}$ as
$\theta^\xi_N:=(i\omega-\cL^\xi_N)^{-1}\gamma_{N,v}^{\xi}$.
Equivalently,
\be\label{cicciobello}
(i\omega-\cL^\xi_N)\theta^\xi_N=\gamma_{N,v}^{\xi}\,.
\en
Notice that we have omitted the dependence on $v$ in $\theta^\xi_N$ to ease the notation.
By writing $\theta^\xi_N(x)= \theta^\xi_{N,R}(x)+ i \theta^\xi_{N,I}(x)$, where $R$ and $I$ refer to real and imaginary parts, we have that \eqref{cicciobello} is equivalent to the system
	\be\label{olleboiccic}
	\begin{cases}
		-\omega \theta^\xi_{N,I}- \cL^\xi _N \theta^\xi_{N,R} =\gamma_{N,v}^{\xi}\,,\\
		\omega \theta^\xi_{N,R}-\cL^\xi_N\theta^\xi_{N,I}=0\,.
	\end{cases}
	\en
Equivalently, $\theta^\xi_{N,R}$ and $ \theta^\xi_{N,I}$ are the unique real functions on $\bbT^d_N$ such that 
 \be\label{pimpa}
 	\begin{cases}
 	-\omega \la f, \theta^\xi_{N,I}\ra_{L^2(\mfm_N)}+\frac 12 \la \nabla f, \nabla   \theta^\xi_{N,R}  \ra_{L^2(\nu_N^\xi)}  
 		= \la f, \gamma_{N,v}^\xi\ra_{L^2(\mfm_N)} \,,\\
 	\omega \la f, \theta^\xi_{N,R}\ra_{L^2(\mfm_N)}+\frac 12 \la \nabla f, \nabla   \theta^\xi_{N,I}  \ra_{L^2(\nu_N^\xi)}  =0
 \end{cases}
\en
for all functions $f:\bbT^d_N \to \bbR$. Indeed, it is enough to test against real functions, but one can take also $f:\bbT^d_N\to \bbC$.

\begin{Lemma}\label{violino} For all $\xi\in \O$ satisfying \eqref{req1}, and in particular for all $\xi\in\O_{\rm typ}$, it holds
  \begin{align}
&  \sup_{N\geq 1}  \|\theta^\xi_N \|_{L^2(\mathfrak{m}_N)} <+\infty\,, \qquad \label{primo}\\
&   \sup_{N\geq 1 } \ \la  \theta^\xi_{N,R} , -\cL^\xi _N \theta^\xi_{N,R} \ra_{L^2(\mathfrak{m}_N) }<+\infty \,, \label{secondo}\\
     &   \sup_{N\geq 1 }  \la  \theta^\xi_{N,I} , -\cL^\xi _N \theta^\xi_{N,I} \ra_{L^2(\mathfrak{m}_N) }<+\infty\,.\label{terzo}
    \end{align}
\end{Lemma}
\begin{proof}
		By taking the scalar product of the left and right hand sides of \eqref{cicciobello} with $\theta^\xi_N$ one obtains
    \begin{align}\label{pouny}
    	i\omega \| \theta^\xi_N\|_{L^2(\mathfrak{m}_N) }^2  +\la  \theta^\xi_N, -\cL^\xi_N\theta^\xi_N\ra_{L^2(\mathfrak{m}_N) }=\la \theta^\xi_N, \gamma_{N,v}^{\xi}\ra_{L^2(\mathfrak{m}_N) }\,.
    \end{align}
    By considering separately the real and imaginary part, the above equation can be rewritten as the system 
    \begin{align}\label{baudo}
    \begin{cases}
    \omega \| \theta^\xi_N\|_{L^2(\mathfrak{m}_N) }^2 
    	=  -\la \theta^\xi_{N,I},\gamma_{N,v}^{\xi}\ra_{L^2(\mathfrak{m}_N) } \,,\\
    \la \theta^\xi_{N,R}, -\cL^\xi_N \theta^\xi_{N,R}\ra_{L^2(\mathfrak{m}_N) }+
    \la \theta^\xi_{N,I}, -\cL^\xi_N  \theta^\xi_{N,I}\ra_{L^2(\mathfrak{m}_N) }
    	= \la \theta^\xi_{N,R},\gamma_{N,v}^{\xi}\ra_{L^2(\mathfrak{m}_N) } \,.
    \end{cases}
    \end{align}
Since $\g_{N,v}^\xi (x) =\sum _{z\in \cN} c_{x,x+z}^{(N)} (\xi) z\cdot v$ and 
$c_{x,x+z}^{(N)}(\xi)=c_{x+z,x}^{(N)}(\xi)$,  for any $h:\bbT^d_N\to \bbC$
by integration by parts  and afterwards by  Cauchy--Schwarz inequality for complex functions we get  \be\label{parti1}
 \begin{split}
  \Big| \langle h,\g_{N,v}^\xi\rangle _{L^2(\mathfrak{m}_N)}  \Big| &=\frac{1}{2N^d} 
  \Big | 
  \sum_{x\in \bbT^d_N} \sum_{z\in \cN} c_{x,x+z}^{(N)}  (\xi) (z\cdot v)  \big( \bar h(x)-\bar h(x+z)\big)
  \Big|\\
 &\leq  \frac{1}{2}\Big( \frac{1}{N^d} \sum_{x\in \bbT^d_N} \sum_{z\in \cN} c_{x,x+z}^{(N)} (\xi) (z\cdot v)^2\Big)^{1/2} \| \nabla h \|_{L^2(\nu^\xi_N)}
 \,.
 \end{split}
 \en

The gradient  $\nabla h$ appearing in \eqref{parti1} is the one  defined in \eqref{grillo}. Since $\xi$ satisfies \eqref{req1}, we have
   \begin{equation*}    \begin{split}
   \frac{1}{N^d} \sum _{x\in \bbT^d_N}  \sum_{z\in \cN} c_{x,x+z}^{(N)}(\xi) 
   \leq   \frac{2}{N^d} \sum _{x\in [0,N)^d\cap \bbZ^d} \sum_{z\in \cN} c_{x,x+z}(\xi) \xrightarrow{N\to\infty} 2\sum_{z\in \cN} \bbE[c_{0,z}]\,,
    \end{split}
    \end{equation*}
    where the factor $2$ is to take care of the border terms.
      Hence, from  \eqref{parti1} we get  for some deterministic constant $C>0$ that 
 \be\label{prop_H_meno}
 \big| \langle h, \g_{N,v}^\xi \rangle _{L^2(\mathfrak{m}_N)}\big|  \leq C \| \nabla h \|_{L^2(\nu_N^\xi)}
 \en
 for any $h:\bbT^d_N\to \bbC$. Then from \eqref{rom} and   \eqref{baudo} we obtain
 \begin{equation}\label{baudobis}
 \begin{cases}
  \o  \| \theta^\xi_N  \|^2_{L^2(\mathfrak{m}_N)} \leq C \| \nabla \theta ^\xi _{N,I}\|_{L^2(\nu^\xi_N)}\,,\\
  \| \nabla \theta^\xi_{N,R}\|^2_{L^2(\nu^\xi_N)} + \| \nabla \theta^\xi_{N,I}\|^2_{L^2(\nu^\xi_N)} \leq 2 C\| \nabla \theta^\xi_{N,R}\|_{L^2(\nu_N^\xi)}\,.
 \end{cases}
 \end{equation}
 The second equation in \eqref{baudobis} implies that   $\| \nabla \theta^\xi_{N,R}\|_{L^2(\nu^\xi_N)} \leq 2 C$ and therefore also $ \| \nabla \theta^\xi_{N,I}\|^2_{L^2(\nu^\xi_N)} \leq 4C^2$.  Using also \eqref{rom}, this last bound and the first equation in \eqref{baudobis} allow to conclude.
 \end{proof}

\begin{Remark}Using that $\bbE[c_{0,z}^2]<+\infty$ for all $z\in \cN$, one can give a more direct proof of Lemma~\ref{violino} by replacing in \eqref{req1} the conductances with the squared conductances. Indeed, from \eqref{baudo} we get
   \begin{equation*}
   \begin{cases}
   & \omega \| \theta^\xi_N\|_{L^2(\mathfrak{m}_N) }^2 
    	\leq  \| \theta^\xi_{N,I}\|_{L^2(\mathfrak{m}_N) } \| \gamma_{N,v}^{\xi}\|_{L^2(\mathfrak{m}_N) }\,,\\
	 & \la \theta^\xi_{N,R}, -\cL^\xi_N \theta^\xi_{N,R}\ra_{L^2(\mathfrak{m}_N) }+
    \la \theta^\xi_{N,I}, -\cL^\xi_N  \theta^\xi_{N,I}\ra_{L^2(\mathfrak{m}_N) }
    \leq  \| \theta^\xi_{N,R}\|_{L^2(\mathfrak{m}_N) } \| \gamma_{N,v}^{\xi}\|_{L^2(\mathfrak{m}_N) }\,. 
    \end{cases}
    \end{equation*}
The first equation implies that  $\| \theta^\xi_N\|_{L^2(\mathfrak{m}_N) }
    \leq \frac{1}{\omega}\| \gamma_{N,v}^{\xi}\|_{L^2(\mathfrak{m}_N) }$ and therefore the r.h.s. in the second equation is upper bounded by  $\frac{1}{\omega}\| \gamma_{N,v}^{\xi}\|_{L^2(\mathfrak{m}_N) }^2$. The claim then follows for all $\xi$ 
   satisfying \eqref{req2} with squared conductances,  since in this case one easily obtains that 
    $\limsup_{N\to+\infty}\| \gamma_{N,v}^{\xi}\|_{L^2(\mathfrak{m}_N) }<+\infty$. \end{Remark}

    Let $\e=1/N$, $N\geq 2$.   We define $\bbT^d_\e$ as the discrete torus  contained in $\bbT^d:=\bbR^d/\bbZ^d$ given by the image of $\bbT^d_N$ under the map  $x\mapsto x/N=\e x$. Hence we have the bijection
    \be\label{bigio}
    \bbT^d_N \ni x \mapsto\e x \in \bbT^d_\e\,.
    \en
    Recall that $\L=(0,1)^d$. The canonical projection $\pi: \bbR^d\to \bbR^d/\bbZ^d=\bbT^d$ is injective when restricted to $\L$.
Recall  that $ \L_\e :=\L \cap ( \e \bbZ^d )$ (cf.~\eqref{mango}). Moreover, recall the definition \eqref{onde_alpha}  of $\mu_\e$ and $\nu_\e^\xi$.
\begin{Definition}\label{pioggia25}
Given $\xi\in \O$ and $\e=1/N$, we define  the functions $\theta^\xi_\e, \psi^\xi_\e , \rho^\xi_\e\in L^2( \mu_\e)$ as 
\[
\theta^\xi_\e ( \e x):= \theta _N^\xi ( \pi_N(x))\,,\qquad 
 \psi^\xi_\e (\e x):= \e \theta^\xi_\e (\e x) \,, \qquad
 \rho^{\xi}_\e( \e x):= \gamma ^{\xi}_{N,v} (\pi_N(x) )\,,
 \]
for all $\e x \in  \L_\e$,  where $\pi_N:\bbZ^d \to \bbZ^d / N \bbZ^d=\bbT^d_N$ is the canonical projection. 
\end{Definition}
Recall the definition of the $\e$--gradient in \eqref{gradientello}.
 As a simple  consequence  of Lemma~\ref{violino}  we have the following:
\begin{Lemma}\label{lemma_stime}
 Given $\z \in \O_{\rm typ}$ the following holds 
  \begin{align}
   & \sup_{\e} \| \theta^{\z}_\e \|_{L^2( \mu_\e)} <+\infty\,, \label{fritz}\\
    & \lim_{\e\da 0} \| \psi^{\z}_\e \|_{L^2(\mu_\e)} =0\,, \label{fritz_bis}\\
& \sup_{\e} \| \nabla_\e \psi^{\z}_\e \| _{L^2 (\nu^\z_\e)} <+\infty\,. \label{kurt}
  \end{align}
\end{Lemma}
\begin{proof}Since $\e=1/N$,   $\e x\in \L_\e=\L \cap \e \bbZ^d$ if and only if $x\in (0,N)^d \cap \bbZ^d$. This implies that $\| \theta^{\z}_\e \|_{L^2(\mu_\e)}\leq \| \theta^{\z}_N \|_{L^2(\mathfrak{m}_N)}$.  By \eqref{primo} in Lemma \ref{violino} we then get \eqref{fritz} and \eqref{fritz_bis}. In addition (cf.~\eqref{onde_alpha}, \eqref{rom} and recall the observation at the beginning of this section) we have 
\begin{equation*}
\begin{split}
\| \nabla_\e \psi^\z_\e \| _{L^2 (\nu^\z_\e)} ^2 & = 
\e^d \sum_{\e x\in \L_\e} \sum_{\substack{z\in \cN:\\ \e x+\e z\in \L_\e }}
c_{x,x+z} (\z)
(\nabla_\e  \psi^\z_\e ) ^2( \e x,z)\\
&= \frac{1}{N^d} \sum_{ x\in (0,N)^d \cap \bbZ^d } \sum_{\substack{z\in \cN:\\  x+ z\in (0,N)^d \cap \bbZ^d  }}c^{(N)}_{ x, x+ z} (\z)
\Big( \theta^\z_N\big( \pi_N(x+ z) \big)-  \theta^\z_N\big( \pi_N( x)\big) \Big)^2\\
& \leq \langle \nabla \theta^\z_N, \nabla \theta^\z_N\rangle _{L^2(\nu^\z_N)} = 2 \langle   \theta^\z_N, -\cL^\z_N \theta^\z_N \rangle_{L^2(\mathfrak{m}_N)}\,.
\end{split}
 \end{equation*}
This estimate together with \eqref{secondo} and \eqref{terzo} (recall \eqref{romix}) implies \eqref{kurt}.
\end{proof}

    As an immediate consequence of \eqref{pimpa}  we have the following result:
    \begin{Corollary}\label{giugno} 
    	Given $\xi \in \O$ and $\d\in (0,1/2)$, for some positive constant $\e(\d)$  it holds
\begin{align}
	-\omega \la f, \theta^{\xi}_{\e,I}\ra _{L^2(\mu_\e)} + \frac{\e^2}{2} \la \nabla_\e f, \nabla_\e \theta^{\xi}_{\e,R} \ra_{L^2(\nu_\e^\xi)}
		&= \la f, \rho_\e^{\xi}\ra_{L^2(\mu_\e)}\label{armando1}\\
	\omega \la f, \theta^{\xi}_{\e,R}\ra _{L^2(\mu_\e)} + \frac{\e^2}{2} \la \nabla_\e f, \nabla_\e \theta^{\xi}_{\e,I}\ra_{L^2(\nu_\e^{\xi})}
		&=0 \label{armando2}
\end{align}
  for any $f:\L_\e\to \bbR$ with support in $\L_\e \cap [\d, 1-\d]^d$ and for all $\e\leq \e(\d)$.
      \end{Corollary}

Furthermore, as a consequence of Lemma  \ref{compatto1}, Lemma \ref{compatto2} and Lemma \ref{lemma_stime}, we also have: 
\begin{Corollary}\label{caciotta}
 Given $\z \in \O_{\rm typ}$,  there exists a subsequence $(\e_k)$ (dependending on $\z$) such that
\begin{align}
& L^2(\L, \mu_\e)\ni \theta^{\z}_\e(x)   \stackrel{2}{\rightharpoonup} \theta^\zeta(x,\xi) \in L^2( \L\times \O, {\rm d}x \times \bbP)\,, \label{salame1}\\
& L^2(\L, \mu_\e)\ni \psi^{\z}_\e(x)   \stackrel{2}{\rightharpoonup} 0 \in L^2( \L\times \O, {\rm d}x\times \bbP)\,, \label{salame2}\\
&  L^2(\L\times \cN, \nu^\xi_\e)\ni \nabla_\e \psi_\e^{\z}(x,z)   \stackrel{2}{\rightharpoonup} w^\zeta(x,\xi,z) \in L^2 ( \L\times \O \times \cN, {\rm d}x \times \nu)\,,  \label{salame3}
\end{align}
for suitable functions $\theta^\zeta, w^\zeta$.
\end{Corollary}

We write $\theta^\z_R,\theta^\z_I$, $w^\z_R$, $w^\z_I$ for the real and the imaginary parts of the functions $\theta^\zeta,w^\zeta$ appearing in \eqref{salame1} and \eqref{salame3}.
 It is convenient to write 
\[ \theta ^\z_x (\xi):=\theta^\zeta (x, \xi)\qquad w^\z_x(\xi,z):= w^\zeta(x,\xi, z) \]
so that when we write $\nabla \theta^\z _x (\xi,z)$ we apply the notation \eqref{cantone} to the function $\xi \mapsto \theta^\z _x(\xi)$. Similarly, when we write $\dive w^\z_x (\xi)$ we apply  the  notation \eqref{emma_0}  to the function $ (\xi,z)\mapsto w^\z_x(\xi,z) $.

Recall that we have fixed a vector $v\in \bbR^d\setminus \{0\}$. Let us consider the form $u_v\in L^2(\O\times \cN,\nu)$ given by
\be\label{tiglio} u_v(\xi,z):= z\cdot v\,.
\en
Note that 
\be\label{algebra}
\frac{1}{2}\dive u_v (\xi) = \g_v(\xi)\,.
\en

The following lemma follows by testing \eqref{armando1} and \eqref{armando2} with test functions of the form $f(x):=\varphi(x)g(\tau_{ x/\varepsilon}\zeta)$
  	with $\varphi\in \mathcal C^\infty_c(\Lambda)$ and $g\in \cG_2$ (in the same spirit  of e.g.~\cite[Claim~19.3]{F_hom} and  \cite{ZP}):
\begin{Lemma}\label{fvg}
For any $\zeta\in\O_{\rm typ}$, $\rmd x$--a.e.~$x\in \L$ and $\bbP$--a.a.~$\xi$ it holds 
\begin{align}
 \o \theta^\z _R(x,\xi)& = \frac{1}{2} \dive w^\z_{x,I}(\xi)\,, \label{fvg1}\\ 
-\o \theta^\z _I(x,\xi)&=\frac{1}{2} \dive \big(w^\z_{x,R}+u_v\big)(\xi)\,.\label{fvg2}
\end{align}

\end{Lemma}
  \begin{proof}
  	We begin with \eqref{fvg1}. 
  	We use  \eqref{armando2} with $\xi$ replaced by $\z$, i.e. 
	\be\label{armando200}
	\omega \la f, \theta^{\zeta}_{\e,R}\ra _{L^2(\mu_\e)} + \frac{\e^2}{2} \la \nabla_\e f, \nabla_\e \theta^{\zeta}_{\e,I}\ra_{L^2(\nu_\e^\zeta)}=0\,,
	\en
	 applied to a test function $f:\,\Lambda_\e\to\bbR$ of the form 
  	$f(x):=\varphi(x)g(\tau_{ x/\varepsilon}\zeta)
  	$
  	with $\varphi\in \mathcal C^\infty_c(\Lambda)$ and $g\in \cG_2$. By the Leibniz rule \eqref{leibniz} the $\e$-gradient of $f$ along $z$  can be written as
  	\begin{align*}
  	\nabla_\varepsilon f(x,z)
  		&=\nabla_\e \varphi(x,z)g(\tau_{x/\e+z}\zeta)+\varphi(x)\big(\nabla_\e g(\tau_{\cdot/\e}\zeta)\big)(x,z)\\
  		&=\nabla_\e \varphi(x,z)g(\tau_{x/\e+z}\zeta)+\varphi(x)\frac{1}{\e}\nabla g(\tau_{x/\e}\zeta,z)
  	\end{align*}
  	where the last gradient of $g$ is in the sense of \eqref{cantone}. Recall that $\psi^\z_\e(\e x)= \e \theta ^\z_\e (\e x)$.
  	Hence we can write the left hand side of  \eqref{armando200} as the sum of the three following terms:
		\begin{align*}
  	A &= \omega \e^d\sum_{x\in\Lambda_\e}\varphi(x)g(\tau_{x/\e}\zeta)\theta^\zeta_{\e,R}(x)\\
  	B &= \frac{\e}{2}\e^d\sum_{ x\in \Lambda_\e}\sum_{\substack{z\in\cN:\\x+ \e z\in \L_\e}}c_{x/\e,x/\e+z}(\z) \nabla_\e\varphi(x,z) g(\tau_{x/\e+z}\zeta) \nabla_\e  \psi_{\e,I}^\zeta(x,z) \\
  	C &= \frac{1}{2}\e^d\sum_{ x\in \Lambda_\e}\sum_{\substack{z\in\cN:\\x+ \e z\in \L_\e}}c_{x/\e,x/\e+z}(\z) \varphi(x)\nabla g(\tau_{x/\e}\zeta,z)\nabla_\e \psi_{\e, I}^\zeta(x,z)\,.
  	\end{align*}
	
  	Taking the limit for $\e\to 0$ along the sequence $(\e_k)$  of Corollary \ref{caciotta} and using  that $\z\in \O_{\rm typ}$ together with   \eqref{rabarbaro} (recall that  $g\in \cG_2\subset \cG$) and the 2-scale convergence \eqref{salame1},  we see that  the term $A$ converges to 
  	\begin{align*}
  	\omega \int_\Omega{\rm d}\bbP(\xi)\int_{\Lambda} {\rm d} x\,\varphi(x)g(\xi)\theta^\z_R(x,\xi)\,.
  	\end{align*}

	Let us now show that the term $B$ goes to $0$ as $\e\to 0$ along  the sequence $(\e_k)$.
 To this aim it is enough to show that 
  \be\label{lavatrice}
 \e^d\sum_{ x\in \Lambda_\e}\sum_{\substack{z\in\cN:\\x+ \e z\in \L_\e}}c_{x/\e,x/\e+z}(\z) \nabla_\e\varphi(x,z) g(\tau_{x/\e+z}\zeta) \nabla_\e  \psi_{\e,I}^\zeta(x,z)
 \en
converges to a finite constant.  
We deal with $B$ using the same arguments  developed in \cite{F_hom} to deal with the first term in \cite[Eq.~(159)]{F_hom}. First,
   we prove that we can replace $\nabla_\e\varphi(x,z)$ with $\nabla \varphi(x)\cdot z$ in \eqref{lavatrice} by showing that the quantity
   \be\label{lavatricebis}
 \e^d\sum_{ x\in \Lambda_\e}\sum_{\substack{z\in\cN:\\x+ \e z\in \L_\e}} c_{x/\e,x/\e+z}(\z) \big(\nabla_\e\varphi(x,z) -\nabla \varphi(x)\cdot z\big) g(\tau_{x/\e+z}\zeta) \nabla_\e  \psi_{\e,I}^\zeta(x,z)
 \en
 goes to zero as $\e \da 0$. 
   We proceed as follows. Being $g\in \cG_2$ a bounded function, by Cauchy-Schwarz inequality we can upper bound \eqref{lavatricebis} by  
   \be\label{lavatricetris}
   \|g\|_\infty \left(\int \rmd\nu_\e^\zeta (u,z) [\nabla_\e \varphi (u,z)-\nabla \varphi (u) \cdot z]^2\right)^{1/2} \|\nabla_\e \psi_{\e,I}^\zeta\|_{L^2(\nu_\e^\z)}\,.
   \en
Since $\z\in \O_{\rm typ}$,   $ \|\nabla_\e \psi_{\e,I}^\zeta\|_{L^2(\nu_\e^\z)}$  is uniformly bounded (cf.~\eqref{kurt}). Applying  Lemma~\ref{pesciolino} now yields that \eqref{lavatricetris} - and therefore also 
\eqref{lavatricebis} -  goes to zero as $\e \da 0$.  
    As a consequence, it remains now to prove  that, for all directions $i=1,2,\dots, d$,  
\be\label{lavatrice4}
 \e^d\sum_{ x\in \Lambda_\e}\sum_{\substack{z\in\cN:\\ x+\e z \in \L_\e}}c_{x/\e,x/\e+z}(\z)\partial_i \varphi(x) z_i  g(\tau_{x/\e+z}\zeta) \nabla_\e  \psi_{\e,I}^\zeta(x,z)
 \en
converges to a finite constant as $\e\to 0 $  along the sequence $(\e_k)$.
	To this aim let $h:\O\times \cN\to \bbR$ be defined as $h(\xi,z):= g(\t_z \xi) z_i$. Note that $h\in \cH_3\subset \cH$ (cf. Section~\ref{sec_typ_env}) and that $h(\t_{x/\e} \z, z)= g( \t_{x/\e +z} \z) z_i$. In particular, \eqref{lavatrice4} equals 
	\be \label{lavatrice5}
	\int_{\L\times \cN} \rmd \nu_\e^\z(x,z) \partial_i \varphi(x)h(\t_{x/\e} \z, z)\nabla_\e  \psi_{\e,I}^\zeta(x,z)\,.
	\en
	By using \eqref{yelena} and  the 2-scale convergence \eqref{salame3}, we conclude that \eqref{lavatrice5} converges as $\e \da 0$ along $(\e_k)$  to $\int_\L \rmd x \int _{\O\times \cN} \rmd \nu(\xi, z) \partial_i \varphi(x)h(\xi,z) w^\z_I(x,\xi, z)$, which is finite since $w^\z\in L^2(\L\times \O\times \cN, \rmd x\times \nu)$ and $h$ is bounded.

Let us now take the limit for $\e\to 0$ along the sequence $(\e_k)$ for the term $C$. Thanks to the 2-scale convergence \eqref{salame3}, using  that $\z\in \O_{\rm typ}$, that  $\nabla g\in  \cH_2\subset \cH$ for all $g\in \cG_2$ and by  \eqref{yelena}, we have that $C$ converges to 
	 \begin{multline*}
     \frac{1}{2} \int_\L {\rm d}x \varphi(x)  \int _{\Omega\times \cN} \rmd\nu (\xi,z) \nabla g(\xi,z)w^\z_{x,I}(\xi,z)
	 =\\ -\frac{1}{2} \int_\L {\rm d}x \varphi(x)  \int _{\Omega} {\rm d}\bbP (\xi)  g(\xi)\dive w^\z _{x,I}(\xi)\,.	
     \end{multline*}
Note that the last identity follows from  \eqref{italia} since $g\in \cG_2$ is measurable and bounded, while for $\rmd x$--a.a. $x\in \L$ we have that $w_{x,I}(\xi,z)\in L^2(\nu)$ since $w(x,\xi, z) \in L^2( \L\times \O\times\cN, \rmd x\times \nu)$.
	
	Let us assemble all the pieces. Since the l.h.s.~of \eqref{armando200} equals $A+B+C$, we get that the limit of $A+B+C$ as $\e \to 0$ along $(\e_k)$ is zero. Hence, we have proved that for any $\varphi\in C_c^\infty(\L)$ and $g\in \cG_2$ it holds
	\be
	\int_{\Lambda}{\rm d} x\, \varphi(x) 
	 \int_\Omega{\rm d}\bbP(\xi)\,g(\xi)\left[ \omega \theta^\z_R(x,\xi)-\frac{1}{2}\dive  w^\z_{x,I}(\xi)\right]=0\,.
	 	\en
  	Due to the density of $C^\infty_c(\L) $ in $L^2_R(\L)$ and the   density of $\cG_2$ in $L^2_R(\bbP)$ (cf.~Section~\ref{sec_typ_env}) we get that $\omega \theta^\z_R(x,\xi)- \frac{1}{2} \dive  w^\z_{x,I}(\xi)=0$ for $\rmd x$-a.e.~$x$  in $\L$ and for $\bbP$--a.a. $\xi$ in $\O$. This proves \eqref{fvg1}.
  	
  	\smallskip
  	
We now prove  \eqref{fvg2}. 
	We use  \eqref{armando1} with $\xi$ replaced by $\z$, i.e. 
	\be \label{armando100}-\omega \la f, \theta^{\z}_{\e,I}\ra _{L^2(\mu_\e)} + \frac{\e^2}{2} \la \nabla_\e f, \nabla_\e \theta^{\z}_{\e,R} \ra_{L^2(\nu_\e^\z)}= \la f, \rho_\e^{\z}\ra_{L^2(\mu_\e)}\,.
	\en
 Using test functions $f$ of the same form as in the previous case, we obtain in a completely similar way that the l.h.s.~of \eqref{armando100} converges  along $(\e_k)$ to
\be\label{abete}
	 \int_{\Lambda}{\rm d} x \, \varphi(x) 
	\int_\Omega{\rm d}\bbP(\xi)g(\xi)\left[ -\omega \theta^\z_I(x,\xi)-\frac{1}{2}\dive  w^\z_{x,R}(\xi)\right]\,.
	\en
	Since  $\varphi \in C^\infty_c(\L)$, the support of $\varphi$ is contained in  the box $ [\d,1-\d]^d$  for  some $\d>0$. As a consequence, for $\e$ small enough,  if $x\in \L_\e$ and $\varphi(x)\not =0$ then  $\gamma ^{\zeta}_{N,v} (\pi_N(x/\e))=\gamma_v (\t_{x/\e} \z)$.  
	Hence, the r.h.s.~of \eqref{armando100} can be rewritten as
	\be \label{bosco}
	\la f, \rho_\e^{\zeta}\ra_{L^2(\mu_\e)}
		=\e^d\sum_{x\in\L_\e}\varphi(x)g(\tau_{x/\e}\zeta)\g_v(\tau_{x/\e}\zeta)\,.
	\en
    By Definition~\ref{def_om_typ} and since $c_{0,z}\in \cG$ for all $z\in \cN$ and $\cG_2\subset\cG$,  we know that  $\O_{\rm typ}\subset   \cA[ g c_{0,z} ] $ for all $z\in \cN$\footnote{We point out that, since $g\in \cG_2$ is bounded, $g c_{0,z}\in L^1(\bbP) $ whenever $c_{0,z}\in L^1(\bbP)$. Hence in the above step we did not really use the assumption that $c_{0,z} \in L^2(\bbP)$. It would have been enough to require  in the definition of $\O$ the inclusion $\O\subset \cap_{z\in \cN}\cap _{g\in \cG_2}\cA[ g c_{0,z}]$.}.  By linearity we then  have $\O_{\rm typ}\subset   \cA[ g \g_v] $  for all $v\in \bbR^d$.
As a consequence and since $\frac{1}{2}\dive u_v  = \g_v$
 (cf.~\eqref{algebra}), \eqref{bosco} converges along $(\e_k)$ to
	\be\label{pino}
	\frac 12\int_{\Lambda}{\rm d} x \varphi(x) 
	\int_\Omega{\rm d}\bbP(\xi)g(\xi)\dive u_v(\xi)\,.
	\en
	Since	\eqref{abete} equals \eqref{pino} for all $\varphi \in C^\infty_c(\L)$ and $g\in \cG_2$, by density we get \eqref{fvg2}.
	  \end{proof}

     \begin{Lemma}\label{venezia}
For any $\zeta\in\O_{\rm typ}$ and for $\nu$--a.a.~$(\xi,z)$ it holds 
	\begin{align}\label{ruperto}
		\int_\L \rmd x\, w^\zeta(x,\xi, z)=\int_\L \rmd x\, \nabla \theta^\z_x (\xi,z) \,.
	\end{align}
\end{Lemma}
    \begin{proof} Recall that the functional set $\cH_*$ introduced in Section~\ref{sec_typ_env} is a countable dense subset of $L^2_R(\nu)$ given by bounded functions.
    Consider a function $\varphi \in\mathcal C_c^\infty(\Lambda)$ and a function $b\in \cH_*$. First of all we observe that for $\e$ small (as understood below)
    \begin{align}\label{labrador}
    A:=\e^d \sum_{ x\in \Lambda_\e} \theta_\e^\zeta(x)\varphi (x)\dive b(\tau_{x/\e}\zeta)
        \end{align}
        equals $B+C$ where
           \begin{align*}
    B&:=-\e^d \sum_{ x\in \Lambda_\e} \sum_{\substack{z\in\cN:\\x+\e z\in \L_\e}} c_{x/\e, x/\e+z}(\z)  \nabla_\e \psi_\e^\zeta(x,z)\varphi (x) b(\tau_{x/\e}\zeta,z)\\
    C&:=-\e^{d+1} \sum_{ x\in \Lambda_\e} \sum_{\substack{z\in\cN:\\x+\e z\in \L_\e}} c_{x/\e, x/\e+z}(\z)  \theta_\e^\zeta(x+\e z)\nabla_\e \varphi (x,z)b(\tau_{x/\e}\zeta,z)\,.
    \end{align*}
    Indeed, by Lemma \ref{cane_ciak} and Remark \ref{passaggio},  $A= - \e \int d \nu _\e^\zeta (x,z) \nabla_\e ( \theta_\e^\zeta \varphi     ) (x,z) b ( \t_{x/\e} \zeta, z)$. 
The claim $A=B+C$ then follows from  \eqref{leibniz} and the identity $\psi_\e^\z(x)
=\e \theta _\e^\z(x)$ for all $x\in \L_\e$.

  Let us now determine the limit  of $A,B,C$ along $(\e_k)$. We recall that $\dive b\in  \cG$ for all $b\in \cH_*$ (see  Definition~\ref{gonars}). Recall that $b$ is bounded as $b\in \cH_*$ and that $\theta^\z_x\in L^2(\bbP)$ for $\rmd x$--a.a. $x\in \L$. Hence by Lemma \ref{faticaccia} (see item (ii) there) we have
  \be\label{italico} 
  \int _\O {\rm d} \bbP(\xi)     \theta^\z(x, \xi) \, \dive b(\xi)= - \int_{\O\times \cN} {\rm d} \nu(\xi, z) \nabla \theta^\z_x (\xi, z) b( \xi, z)  \qquad \rmd x\text{--a.s.} \,.
  \en
  Note that to derive \eqref{italico}  we used the assumption $c_{0,z}\in L^2(\bbP)$. 
  Going back to the expression in \eqref{labrador} for $A$, we see by \eqref{salame1} (and therefore \eqref{rabarbaro}) and afterwards by \eqref{italico}   that along the sequence $(\e_k)$ it holds 
    \begin{align*}
    \lim_{\e\da 0}A
    	&=\int_{\Omega}\,\rmd\bbP(\xi)\int_{\Lambda}\rmd x\theta^\zeta(x,\xi)\varphi (x) \dive b(\xi)\\
    	&=- \int_{\Omega\times \cN}\,\rmd\nu(\xi,z)\int_{\Lambda}\rmd x\nabla\theta^\z_x(\xi,z)\varphi (x) b(\xi,z)\,.
    \end{align*}
    On the other hand, \eqref{salame3} and the property $b\in \cH_*\subset \cH$ (cf.~\eqref{yelena}) yields
    \begin{align*}
    \lim_{\e\da 0}B
    	=- \int_{\Omega\times \cN}\,\rmd\nu(\xi,z)\int_{\Lambda}\rmd x \,w^\zeta(x,\xi,z)\varphi (x) b(\xi,z)\,.
    \end{align*}
    
  We finally show that the quantity $C$ vanishes as $\e\to0$. 
  Indeed, we will prove that $C/\e$ converges to a finite quantity.
  By applying \eqref{micio3} in Lemma~\ref{gattonaZ} and Remark \ref{passaggio} we rewrite $C$ as 
   \be\label{sale} C=\e^{d+1} \sum_{ x\in \Lambda_\e} \sum_{\substack{z\in\cN:\\ x+\e z\in \L_\e }} c_{x/\e, x/\e+z}(\z) \theta_\e^\zeta(x)\nabla_\e \varphi (x,z) \tilde b(\tau_{x/\e}\zeta,z)
  \en
  where $\tilde b (\xi,z):= b(\t_z \xi, -z)$. Recall that  $\tilde b \in \cH$ for any $b\in \cH$. 
  We claim that 
  \be\label{usignolo5}
  \lim _{\e \da 0}\e^{d} \sum_{ x\in \Lambda_\e}\sum_{\substack{z\in\cN:\\ x+\e z\in \L_\e }} c_{x/\e, x/\e+z} (\z) \theta_\e^\zeta(x)[\nabla_\e \varphi (x,z)- \nabla \varphi(x) \cdot z]   \tilde b(\tau_{x/\e}\zeta,z)=0\,.
  \en
  Recall that $b$ (and therefore $\tilde b$) is uniformly bounded as $b\in \cH_*$. Moreover, since $\varphi\in C_c^\infty(\L)$, by Taylor expansion $ |\nabla_\e \varphi (x,z)- \nabla \varphi(x) \cdot z| \leq C \e $ uniformly in $x$ ad $z$ as in \eqref{usignolo5}. Hence we just need to show that the expression
 \[ \e^{d} \sum_{ x\in \Lambda_\e}\sum_{\substack{z\in\cN:\\ x+\e z\in \L_\e }}c_{x/\e, x/\e+z}(\z) |\theta_\e^\zeta(x)| 
 \]
 remains bounded as $\e \da 0$. By  Cauchy-Schwarz inequality the above expression is upper bounded by 
 \[
 \sqrt{|\cN|} \| \theta_\e^\zeta  \|_{L^2(\mu_\e)}\Big(  \e^{d} \sum_{ x\in \Lambda_\e}\sum_{\substack{z\in\cN:\\ x+\e z\in \L_\e }}  c^2_{x/\e, x/\e+z}(\z) \Big)^{1/2}
 \]
 and the latter remains bounded as $\e\da 0$ by \eqref{req2} and \eqref{fritz}. Note that in the last argument we have used that $c_{0,z}\in L^2(\bbP)$.

By \eqref{usignolo5}   we can replace in \eqref{sale} $\nabla_\e \varphi(x,z)$ by $\nabla \varphi(x)\cdot z$ when taking $\e \da 0$. Hence, letting
   $h_i(\xi):=\sum_{z\in \cN} c_{0,z}(\xi)  z_i \tilde b(\xi,z) $, we have (using that $\nabla\varphi$ has support in $[\d,1-\d]^d$ for some $\d>0$ and that $c_{0,z} (\t_{x/\e} \z) =c_{x/\e, x/\e+z}(\z)$)
   \[
\frac 1 \e C    = \sum_{i=1}^d \e^{d} \sum_{ x\in \Lambda_\e}  \theta_\e^\zeta(x)\partial_i  \varphi (x)h_i (\tau_{x/\e}\zeta)+o(1)\,.
   \]
   At this point by \eqref{salame1}  and since $h_i \in \cG$ for any $i$ due to Definition~\ref{gonars}, by taking the limit along $(\e_k)$, we get
    \begin{align*}
    \lim_{\e\da 0} \frac 1 \e C
    	&=\sum_{i=1}^d \int_{\Omega}\rmd\bbP(\xi)\int_{\Lambda}\rmd x\,\theta^\zeta(x,\xi) \partial_i \varphi(x) h_i(\xi)    	<\infty\,.
    \end{align*}
 
    To recap, by taking the limit for $\varepsilon\to 0$  along $(\e_k)$ on both sides of the equation $A=B+C$, we have found that
    \begin{align*}
    \int_{\Omega\times \cN}\,\rmd\nu(\xi,z)\int_{\Lambda}\rmd x\nabla\theta^\z_x(\xi,z)\varphi (x) b(\xi,z)
    	=\int_{\Omega\times \cN}\,\rmd \nu(\xi,z)\int_{\Lambda}\rmd x \,w^\zeta(x,\xi,z)\varphi (x) b(\xi,z)\,.
    \end{align*}
    The relation \eqref{ruperto} follows then by the density of $C^\infty_c(\L)$ in $L^2_R(\L)$ and the density of $\cH_*$ in $L^2_R(\nu)$.    \end{proof}

We define 
\be\label{grande_theta} \Theta^\zeta(\xi):=\int _\L \rmd x\, \theta^\zeta (x, \xi)
\en as a function in $L^2(\bbP)$.  Note that the above integrand is integrable for $\bbP$--a.a.~$\xi$  and  the resulting function $ \Theta^\zeta$ belongs to $L^2(\bbP)$. Indeed,   by Cauchy-Schwarz inequality, we have 
\begin{align*}
    & \int_\O\rmd \bbP(\xi) \int _\L \rmd x\, |\theta^\zeta (x, \xi)|\leq \Big(\int_\O\rmd \bbP(\xi) \int _\L \rmd x\, |\theta^\zeta (x, \xi)|^2 \Big)^{1/2}<+\infty\,,\\
    &\int_\O \rmd \bbP(\xi)| \Theta^\zeta(\xi)|^2 
\leq 
\int _\O \rmd \bbP(\xi) \int _\L \rmd x |\theta^\z(x,\xi)|^2<+\infty\,.
\end{align*}

The boundedness in the righ hand sides follows from the fact that $\theta^\z\in L^2( \L\times\O, \rmd x \times \bbP)$. 

We write $\Theta_R^\zeta$ and $\Theta_I^\zeta$ for the real and imaginary parts of $\Theta^\zeta$. 
 For the next result recall that $u_v(\xi,z):=z \cdot v$ (see \eqref{tiglio}):

\begin{Lemma}\label{biauzzo}
 For any $\zeta\in\O_{\rm typ}$ and for $\bbP$--a.a.~$\xi$  it holds
  \begin{align}
    &    \o \Theta^\z_R  (\xi)= \frac{1}{2}\dive \big(\nabla \Theta^\z_I\big)(\xi)  \,,\label{ilpatriota}\\
   &  \o \Theta^\z_I  (\xi)=- \frac{1}{2}  \dive \big(\nabla  \Theta^\z_R +u_v\big) (\xi)
    \label{ilpatriota2lavendetta}\,.
  \end{align}
\end{Lemma}
Above, the gradients $\nabla \Theta^\z_I $ and $\nabla \Theta^\z_R$ are as in \eqref{cantone}.
\begin{proof}
For  \eqref{ilpatriota} we see that $\bbP$--a.s.
\[ 
2 \o \Theta^\z_R   =
	\int _\L \rmd x\, \big(\dive w^\z_{x,I}\big)  = \dive \Big( \int _\L \rmd x\,  w^\z_{x,I} \Big) \\
     = \dive \Big( \int _\L \rmd x\,  \nabla\theta^\z_{x,I}\Big) 
	=\dive (\nabla\Theta^\z_I ) \,,
\]
where the first equality comes from equation \eqref{fvg1} in Lemma \ref{fvg} and the third equality from Lemma \ref{venezia}.

Analogously for \eqref{ilpatriota2lavendetta}, using equation \eqref{fvg2} in Lemma \ref{fvg} and Lemma \ref{venezia}, 
 we have  $\bbP$--a.s.
\begin{align*}
-2\o \Theta^\z_I 
	&=  \int _\L \rmd x\, \dive \big(w^\z_{x,R} +  u_v\big) = \dive\Big(\int _\L \nabla\theta^\z_{x,R}\,\rmd x\Big)+\dive u_v
	= \dive \big(\nabla  \Theta^\z_R + u_v\big) \,.
\end{align*}
\end{proof}

 In the rest $\bbL$ will denote the complex extension of the operator $\bbL$ introduced in Proposition~\ref{zap}.
In the following lemma we  investigate the equations appearing in  Lemma~\ref{biauzzo}.
\begin{Lemma}\label{grovis}  Suppose a function $h \in L^2 (\bbP)$ satisfies
\be\label{delfino98}
  \begin{cases}
    \o h_R =\frac{1}{2} \dive  (\nabla h_I)\,,\\
     \o h_I=-\frac{1}{2} \dive (\nabla h_R+ u_v)\,.
  \end{cases}
  \en
Then $h\in \cD(\bbL)$ and 
  $(i \o - \bbL) h=\gamma_v$.
  \end{Lemma}
\begin{proof}   
The proof is based on the  criterion given at the end of Proposition \ref{zap}  applied to  $h_R$ and $h_I$. 

Recall the operator $L f$ introduced there. Since for a generic function $f=f(\xi)$ and $z\in \cN$ it holds $ \nabla f(\xi,z)- \nabla f(\t_z \xi, -z)= 2 \big( f(\t_z\xi) - f(\xi)\big)$, 
we get  that 
$\frac{1}{2} \dive  (\nabla h_R)= L h_R $ and $ \frac{1}{2} \dive  (\nabla h_I)= L h_I $.
Since, in addition, $\frac{1}{2}\dive u_v (\xi) = \g_v(\xi)$ (cf.~\eqref{algebra}), the system \eqref{delfino98} can be rewritten as
\be\label{balena}
\begin{cases}
    L h_I= \o h_R \,,\\
     L h_R=-\o h_I-\g_v\,.
  \end{cases}
  \en
By \eqref{balena}  and since $h, \g_v \in L^2(\bbP)$ (for $\g_v$  we use $c_{0,z}\in L^2(\bbP)$ for all $z\in \cN$), we get that  $L h_I $ and $L h_R$ belong to $L^2_R(\bbP)$.  This shows that both $h_I$ and $h_R$ satisfy item (ii) in Proposition~\ref{zap}.  They also satisfy item (i) since $c_{0,z}$, $h_I$, $h_R$, $h_I\circ \t_z $, $h_R\circ \t_z$ belong to $L^2_R(\bbP)$ (for the last two functions use the stationarity of $\bbP$). Hence, by Proposition~\ref{zap} we get that $h_I,h_R\in \cD(\bbL)$  and $\bbL h_I= \o h_R $, 
     $\bbL h_R=-\o h_I-\g_v$. This implies  that $h\in \cD(\bbL)$ and $$(i \o - \bbL) h=(i \o - \bbL) (h_R+i h_I)= -\bbL h_R -\o h_I   + i ( \o h_R-\bbL h_I)=\g_v\,.$$
\end{proof}


Let us write $\Theta_k^\zeta$ for the function $\Theta^\zeta$ in \eqref{grande_theta} when $v=e_k$. Similarly we write $\g_k$ for  the function $\g_v$ with $v=e_k$. Trivially, $\g_k$ is the $k$--th component of the function $\g$ appearing in Theorem~\ref{teo1bis} (see \eqref{def_gammaV}).
By combining Lemma \ref{biauzzo} and Lemma \ref{grovis} we have the following:
\begin{Corollary}\label{gioioso} Given $\z\in \O_{\rm typ}$ and $k=1,\dots, d$, the function $\Theta_k^\zeta$ belongs to $\cD(\bbL)$ and satisfies
$(i\o -\mathbb L ) \Theta_k^\zeta=\g_k$.
\end{Corollary}

\begin{Remark}\label{remarkone} For possible further progress  towards weakening the moment assumption, we point out that in this section  we used that $c_{0,z}\in L^2(\bbP)$ only  in the derivation of \eqref{italico} and  \eqref{usignolo5} and in the application of Proposition~\ref{zap} inside the proof of Lemma~\ref{grovis}. Apart from these applications, for all the other arguments in this section we could have avoided to impose $c_{0,z}\in \cG$ by requiring  in Definition~\ref{def_om_typ}  that $\O$ is included in $\cap _{z\in \cN}\cap _{g\in \cG_2}  \cA[c_{0,z} g]$ (imposing that $1\in \cG_2)$.
 \end{Remark}

\section{Proof of Theorem~\ref{teo1bis} (II): Conclusion}\label{sec_conclusione}

In this section, for simplicity of notation, 
we abbreviate $\Lambda_N:=[0,N)^d \cap \bbZ^d$ and, given $f:\bbT^d_N\to \bbC$, we denote again by $f$ the lifting  of $f$ to all of $\bbZ^d$, i.e. $f(x):=f(\pi_N(x))$ for all $x\in \bbZ^d$.

Having the homogenization results of the previous section, we can finally  conclude the  proof of Theorem~\ref{teo1bis}. \rosso{Recall that, by Remark~\ref{nevischio} adapted to Theorem~\ref{teo1bis}, we dot not need to derive  the  second identity in \eqref{jabba2bisV}}.
By comparing \eqref{jabba2*} and \eqref{jabba2bisV} and using Corollary~\ref{gioioso}, to prove Theorem~\ref{teo1bis} it is enough to show that, for all $z\in \cN$, for all $j,k\in\{1,\dots,d\}$  and for all $\z\in \O_{\rm typ}$, 
\begin{align}
&          \lim_{N\to+\infty} \mfm_N[c_{\cdot,\cdot+z}^{(N)}]=   \bbE[c_{0,z}]   \,,\label{limite1} \\
&        \lim_{N\to +\infty} \la \gamma_{N,j}^{\zeta}, (i \o -\cL^{\zeta}_{N})^{-1} \gamma _{N,k}^{\zeta} \ra_{L^2(\mfm_N)}
	=\int _\O \rmd\bbP( \xi) \g_j( \xi)\Theta_k^\zeta ( \xi)         \,.\label{limite2}
\end{align}


Limit \eqref{limite1} is simple. Indeed, by Definition~\ref{def_per_rates}, we have 
\[ \mfm_N[c_{\cdot,\cdot+z}^{(N)}]=\frac{1}{N^d} \sum_{x\in \L_N}  c_{x,x+\bar z}(\z)\]
where $\bar z:= z$ if $z\in \cN_*$ and  $\bar z:=-z$ if $z\in -\cN_*$.
Therefore \eqref{limite1} follows from \eqref{req1} \rosso{and the identity $\bbE[c_{0,z}]=\bbE[c_{0,-z}]$}.

\smallskip

The rest of this section is devoted to the proof of \eqref{limite2}. Since the homogenization results hold along sequences we argue as follows. To prove \eqref{limite2} 
it is enough to show that, for $\zeta\in \O_{\rm typ}$ fixed,  for any diverging sequence $(N_k)_{k\geq 1}$ we can extract a subsequence $(N_{k_r})_{r\geq 1}$  such that  the limit \eqref{limite2} holds  with $N$ in the subsequence $(N_{k_r})_{r\geq 1}$. We notice that, as \rosso{used also} in \cite{F_hom},  the conclusions of Lemma~\ref{compatto1} and Lemma~\ref{compatto2} still hold if one considers at the beginning  $\e$  as vanishing along a sequence (the final convergence in the above lemmas would  then  be along   a subsequence). As a consequence the same holds for \eqref{salame1}, \eqref{salame2} and \eqref{salame3} appearing in Corollary \ref{giugno}, that is, \eqref{salame1}, \eqref{salame2} and \eqref{salame3}  are valid  for $\e$ in a suitable subsequence of a given vanishing sequence $(\e_k)$. 
  
Hence, to simplify the notation we just prove that \eqref{limite2} holds along a diverging sequence (possibly depending on $\zeta$), but all can be restated for a subsequence of a given diverging sequence due to the  comments above on \eqref{salame1}, \eqref{salame2} and \eqref{salame3}.

We first state and prove an auxiliary result \rosso{where $\theta^\z_{N,k}$ is defined as the function $\theta_N^\z$ introduced in Section~\ref{sec_dim_I} with $v=e_k$ (cf.~\eqref{cicciobello}).}
\begin{Lemma}\label{pirano} There exists $C>0$ with the following property.
Let $\psi\in C_c(\bbR)$ be with values in $[0,1]$ and with support in  $[1/2-a,1/2+a]  $  for some $a\in (0,1/2)$.
  Let   $\zeta \in \O_{\rm typ} $, $ k\in \{1,2,\dots, d\}$ and $z\in \cN $. Then
\[ \varlimsup _{N\to +\infty}\frac{1}{N^d} \sum_{x \in \L_N}   \psi\big(\big |\frac{x}{N}-p\big |_\infty\big) c^{(N)}_{x,x+ z}(\z) \big|  \theta^\zeta_{N,k} (x) \big|\\
 \leq C \left( (\tfrac{1}{2}+a)^d - (\tfrac{1}{2}-a)^d\right) 
 \]
where $p$ denotes the center of the box $[0,1]^d$, i.e.\,$\, p:=(1/2,1/2,\dots, 1/2)$.
\end{Lemma}
\begin{proof}
 By Cauchy--Schwarz inequality we can bound 
\[ N^{-d} \sum_{x \in \Lambda_N} \psi\big(\big |\frac{x}{N}-p\big |_\infty\big) c^{(N)}_{x,x+z}(\z)\big| \theta^\zeta_{N,k} (x) \big|\leq D_N (\zeta)^{1/2} E_N(\zeta)^{1/2}
\]
where
\begin{align}
&  D_N(\zeta):=
   N^{-d}  \sum_{x\in \Lambda_N } 
    \psi(|\frac{x}{N}-p|_\infty)^2
   c^{(N)}_{x,x+z}(\zeta) ^2\,, \label{tonkatsu} \\
  & E_N(\zeta):= N^{-d} \sum _{x\in \Lambda_N } \theta^\zeta _{N,k} (x)^2= \| \theta_{N,k}^\zeta \|^2 _{L^2 (\mathfrak{m}_N)}\,.
\end{align}
By \eqref{primo}  and since $\zeta\in \O_{\rm typ}$,  $E_N(\zeta)$ is bounded by a constant $C^2$  uniformly in $N$.

We only need to bound $\varlimsup_{N\to +\infty} D_N(\zeta)$. We treat separately the cases $z\in\cN_*$ and $z\not\in\cN_*$. If $z\in\cN_*$ we have that $c_{x,x+z}^{(N)}(\zeta)=c_{x,x+z}(\zeta)$ so that easily
\begin{align}\label{torello}
D_N(\zeta)
	\leq N^{-d}  \sum_{x\in \bbZ^d }\psi(|\frac{x}{N}-p|_\infty)^2 c^{\,2}_{x,x+z}(\zeta) \,.
\end{align}
Since 
$\zeta \in \O_{\rm typ}\subset  
\cap _{z\in\cN}  \cA [c_{0,z}^2]$ (see Definition~\ref{def_om_typ}), by \rosso{Proposition~\ref{prop_ergodico}} we have 
\be
\lim_{N\to +\infty} N^{-d}  \sum_{x\in \L_N }  \psi(|\frac{x}{N}-p|_\infty)^2
c^{\,2}_{x,x+z}(\zeta) = \int _{\bbR^d} \rmd u \,\psi ( |u-p|_\infty)^2 \bbE[c_{0,z}^2]\,.
\en
To conclude it is enough to observe that  $\int _{\bbR^d} \rmd u\, \psi ( |u-p|_\infty) ^2$ is bounded by $(1/2+a)^d - (1/2-a)^d $.

If $z\not\in\cN_*$ we can split the sum in \eqref{tonkatsu} as two sums over the sets $A_N=\{x\in\L_N:\,x+z\in\Lambda_N\}$ and $B_N=\{x\in\L_N:\,x+z\not\in\Lambda_N\}$. For the sum over $A_N$, we can bound it as in \eqref{torello}  since for $x\in A_N$ we have $c^{(N)}_{x,x+z}(\z)=c_{x,x+z}(\z)$. We move to $B_N$. If $x\in B_N$ then $c^{(N)}_{x,x+z}(\z)= c _{u, u-z}(\z)$, where $u$ is the only element of $\L_N$ such that $\pi_N(x+z)=\pi_N(u)$ (see Definition~\ref{def_per_rates}). Note that $u-z\not \in \L_N$. Hence the contribution to $D_N(\z)$ coming from  the vertexes $x\in B_N$ is upper bounded by
\be\label{temporale}
\frac{1}{N^d} \sum _{u\in \L_N\,:\, u-z\not \in \L_N} {c}^{\,2}_{u,u-z}(\z)
\en
(above we have used that the map $\{x,x+z\}\mapsto \{u,u-z\}$ is injective).
Using that $\z\in \O_{\rm typ}\subset \cA[c^2_{0,z}]$ and applying \rosso{Propositon~\ref{prop_ergodico}} with test functions $\varphi$ whose support is  more and more concentrated on the border of $[0,1]^d$,
 we conclude that  \eqref{temporale}, and therefore also the contribution of $B_N$ to $D_N(\z)$, vanishes as $N\to +\infty$.
\end{proof}

We can now prove the limit \eqref{limite2} along a diverging sequence $(N_k)$.
To this aim  we fix $a\in(0,1/2)$ and take  a continuous function $f:\bbR_+\to [0,1]$  
such that $f(s)=1$ for $s\in [0,1/2-a]$ and $f(s)= 0$ for $s>1/2$. We set $\varphi (u):=f(|u-p|_\infty)$ \rosso{where $p$ is as in Lemma~\ref{pirano}}. Note that  $\varphi(u)=1$ if $u \in [a,1-a]^d$ and $\varphi(u)=0$ if $u\not\in [0,1]^d$.
\rosso{By \eqref{cicciobello} we} can write
\be\label{rondine}
\begin{split}
  &\la \gamma_{N,j}^\zeta, (i \o -\cL^{\zeta}_N)^{-1} \gamma _{N,k}^{\zeta} \ra_{L^2(\mfm_N)} = 
    \la \gamma_{N,j}^{\zeta}, \theta^\zeta_{N,k} \ra_{L^2(\mfm_N)}=N^{-d} \sum_{x\in \L_N} \gamma_{N,j}^{\zeta}(x) \theta^\zeta_{N,k} (x)\\
     & = N^{-d} \sum_{x\in \L_N } \varphi (x/N) \gamma_{N,j}^{\zeta}(x) \theta^\zeta_{N,k} (x)+N^{-d} \sum_{x\in \L_N}\left( 1-\varphi(x/N)\right)  \gamma_{N,j}^{\zeta}(x) \theta^\zeta_{N,k} (x)\\
     &= :A_N(\zeta)+B_N(\zeta)\,.
\end{split}
\en
Note that  (see Definition~\ref{pioggia25})
\begin{align}
A_N(\zeta)
	& =\e^d \sum_{x\in \L_\e}\varphi (x) \sum_{z\in \cN_*}\left[c _{x/\e,x/\e+z}(\zeta)- c _{x/\e,x/\e-z}(\zeta)\right] z_j \theta ^\zeta _{\e,k}(x)\,,\\
|B_N(\zeta) | 
	& \leq C  N^{-d} \sum_{x\in \L_N} (1-\varphi (x/N)) \Big[ \sum_{z\in \cN} c^{(N)}_{x,x+z}(\zeta) \Big] \left| \theta^\zeta_{N,k} (x)\right|\,.
\end{align}
By Lemma \ref{pirano} and since $\zeta\in \O_{\rm typ}$ we have  
\be\label{rondine0}
\varlimsup_{N\to \infty} |B_N(\zeta)|\leq C\rosso{|\cN|}\left( (\tfrac{1}{2}+a)^d-(\tfrac{1}{2}-a)^d\right)\,.
\en

Since $\zeta\in \O_{\rm typ}$ due to \eqref{salame1} and since  $c_{0,z}(\zeta)$ belongs to $\cG$ for each $z\in \cN$ and $c _{x/\e,x/\e\pm z}(\zeta)= c _{0,\pm z}(\t_{x/\e} \zeta)$, we have along a diverging  sequence $(N_k)$ associated to $\zeta$ that 
\be\label{pulcino}
\lim_{N\to +\infty} A_N(\zeta) =
\int _\L \rmd x\,\varphi(x) \int _\O \rmd\bbP( \xi) \g_j( \xi)\theta_k^{\z} (x, \xi)
\en
where $\theta_k^\z(x,\xi)$ is the function in \eqref{salame1} when $v:=e_k$.

Note that, by Cauchy--Schwarz inequality, 
\be\label{rondine1}
\begin{split}
& \left | \int _\L \rmd u \left(1-\varphi(u)\right) \int _\O \rmd\bbP( \xi) \g_j( \xi)\theta^\z_k (x, \xi)
\right|\\
& \leq \left( \int _\L \rmd u (1-\varphi(u) )^2  \right)^{1/2}  \left( \int_\L \rmd u \left(  \int _\O \rmd\bbP( \xi) \g_j( \xi)\theta_k^\z (u, \xi)\right )^2 \right)^{1/2}\,.
\end{split} 
\en
By applying again Cauchy--Schwarz inequality we have 
\be\label{rondine2}
\begin{split}
  \int_\L \rmd u \left(  \int _\O \rmd\bbP( \xi) \g_j( \xi)\theta_k^\z (u, \xi)\right )^2 
 &  \leq 
   \int_\L  \rmd u  \int _\O \rmd\bbP( \xi) \g_j( \xi)^2 \int_\O \rmd\bbP( \xi )  \theta^\z_k (u, \xi)^2\\
   & = \| \g_j\|_{L^2(\bbP)} ^2 \| \theta^\z_k \| _{L^2( \L\times \O, dx \times \bbP)}^2<+\infty\,.
   \end{split} 
\en
By combining \eqref{rondine1} and \eqref{rondine2} we have 
\begin{multline}\label{rondine3}
\left | \int _\L \rmd u\,\varphi(u) \int _\O \rmd \bbP( \xi) \g_j( \xi)\theta^\z_k (u, \xi)-
\int _\L \rmd u  \int _\O \rmd \bbP( \xi) \g_j( \xi)\theta^\z_k (u, \xi)\right|\\
\leq  
\rosso{\left(1- (1-a)^d\right)^{1/2} \| \g_j\|_{L^2(\bbP)}  \| \theta^\z_k \| _{L^2( \L\times \O, dx \times \bbP)} }\,.
\end{multline}
Due to \eqref{rondine},  \eqref{rondine0}, \rosso{\eqref{pulcino}} and \eqref{rondine3} we get
\[
\varlimsup _{N\to +\infty} \left | 
\la \gamma_{N,j}^{\zeta}, (i \o -\cL^{\zeta}_{N})^{-1} \gamma _{N,k}^{\zeta} \ra_{L^2(\mfm_N)}-\int _\L \rmd u  \int _\O \rmd\bbP( \xi) \g_j( \xi)\theta^\z_k (u, \xi)
\right |
 \leq C (a)\,,
 \]
where $a=a(\varphi)$ and  $\lim_{a\da 0} C(a)=0$.

Since we can exhibit $\varphi$ with $a(\varphi)$ arbitrarily small,  we get 
\be\label{eccolo}
\lim_{N\to +\infty} \la \gamma_{N,j}^{\zeta}, (i \o -\cL^{\zeta}_{N})^{-1} \gamma _{N,k}^{\zeta} \ra_{L^2(\mfm_N)}
	=\int _\L \rmd u  \int _\O \rmd\bbP( \xi) \g_j( \xi)\theta^\z_k (u, \xi)\,.
\en
At this point \eqref{limite2} follows from \eqref{eccolo} and definition \eqref{grande_theta}.

\appendix

\section{Proof of Proposition~\ref{zap}}\label{proof_prop_zap}
{\bf Warning:} \emph{To ease the notation, in this appendix  all functional spaces are made by real fuctions. Hence we write i.e. $L^2(\bbP)$ instead of $L^2_R(\bbP)$.}

 Since $\bbE[c_{0,z}(\xi)]<+\infty$ for all $z\in\cN$, the domain $\cD(\cE)$ contains  the subspace of bounded functions, which is dense in $L^2(\bbP)$. Using also the definition \eqref{def_dir_form}, we get that  $\cE$ is a symmetric form on $L^2(\bbP)$. Let us prove that it is closed. To this aim let  $(f_n )$ be a sequence in $\cD(\cE)$ which is Cauchy for the norm $\|\cdot \|_1$ defined before Proposition~\ref{zap}. Being Cauchy in $L^2(\bbP)$, at the cost of extracting a subsequence, we get that $f_n \to f $ in $L^2(\bbP)$ and $\bbP$--a.s., for some function $f\in L^2(\bbP)$.
By Fatou's lemma and due to the $\bbP$--a.s. convergence we get 
$ \cE(f- f_m,f-f_m) \leq\varliminf_{n\to \infty} 
 \cE(f_n- f_m,f_n-f_m)$\,. Hence, since $(f_n)$ is Cauchy w.r.t. $\|\cdot\|_1$, we conclude that $\lim_{m\to \infty} \cE(f- f_m,f-f_m) =0$. 
The above observation  and the fact that $f_m\to f$ in $L^2(\bbP)$ imply that 
  $\lim_{m\to \infty} \| f-f_m\|_1=0$, thus concluding  the proof that $\cE$ is closed. To prove that $\cE$ is a Dirichlet form we need to show that it is Markovian (see the discussion before Proposition~\ref{zap}). This is an immediate consequence of the definition of normal contraction (see \eqref{zompicchia} and of \eqref{def_dir_form}).
  
  \smallskip
  
  Let us now prove that   the set $\cV$ of measurable bounded   functions 
 is a core for $\cE$. Trivially, $\cV$   is included in $\cD(\cE)$ as $c_{0,z}\in L^1(\bbP)$ for all $z\in \cN$. We need to show that for any $f\in \cD(\cE)$ there exists a sequence $(f_n)$ in $\cV$ such that  $\|f-f_n\|_{L^2(\bbP)} \to 0 $   and $\cE(f-f_n,f-f_n)\to 0$ as $n\to +\infty$. To this aim it is enough to take $f_n$ as the $n$--cutoff function 
 \be\label{ncutoff}
 f_n(\xi):= \mathds{1}(|f(\xi) |\leq n) f(\xi) -n \mathds{1}( f(\xi) <-n) +n \mathds{1} (f(\xi) >n)\,.
 \en 
 By dominated convergence  we get  $\|f-f_n\|_{L^2(\bbP)} \to 0 $. To prove that  $\cE(f-f_n,f-f_n)\to 0$, we observe that $|f_n(\xi)-f_n(\xi')| \leq |f(\xi)-f(\xi')|$.
 It then follows that 
$  \left( (f-f_n) (\t_z \xi) - (f-f_n)(\xi)\right)^2    \leq 4  \left( f (\t_z \xi) - f(\xi)\right)^2 $. The above bound and the fact that $f\in \cD(\cE)$ allow to use 
 dominated convergence and  to   conclude that $\cE(f-f_n,f-f_n)\to 0$, since $f-f_n \to 0$ pointwise as $n\to +\infty$.

\smallskip

 Since $\cE$ is a closed symmetric form,  existence and uniqueness 
  of a negative semidefinite self-adjoint operator $\bbL:  \cD(\bbL)\to L^2(\bbP)$ with $\cD(\bbL)\subset L^2(\bbP)$  and satisfying \eqref{istria} follow from \cite[Theorem~1.3.1]{FOT}.

\smallskip

Let us now prove the criterion in Proposition~\ref{zap} assuring that $f\in \cD(\bbL)$ and $\bbL f = Lf $. Take $f\in L^2(\bbP)$ satisfying items (i) and (ii). Let us call $\tilde{\cE}(f,g)$ the r.h.s.~of \eqref{def_dir_form} when  the function inside the expectation belongs to $L^1(\bbP)$.
 The proof of the criterion will consist in the following four steps: 
 
 (a) we prove that $\tilde\cE(f,g)= - \la Lf , g\ra_{L^2(\bbP)}$ for all $g\in \cV$, 
 
 (b) we use the previous property to deduce that $f\in \cD(\cE)$, 
 
 (c) by the previous two properties we show that $\cE(f,g)= - \la Lf , g\ra_{L^2(\bbP)}$ for all $g\in \cD(\cE)$, 
 
 (d) we conclude by applying \cite[Proposition 2.12]{R}.

\smallskip

Let us start.
 We take $g\in \cV$. 
 By item (i), the maps $\xi \mapsto c_{0,z}(\xi) \left( f (\t_z \xi) - f(\xi)\right)  g (\t_z \xi) $ 
 and $\xi \mapsto  c_{0,z}(\xi) \left( f (\t_z \xi) - f(\xi)\right)  g(\xi)$ are both in $L^1(\bbP)$. 
Moreover, using the stationarity of $\bbP$ and that $c_{0,z}(\t_{-z}\xi)= c_{-z,0}(\xi)=c_{0,-z}(\xi)$ we have 
\be\label{NY}
\bbE\big[ c_{0,z}(\xi) \left( f (\t_z \xi) - f(\xi)\right)  g (\t_z \xi)\big] =-\bbE\big[ c_{0,-z}(\xi) \left( f (\t_{-z} \xi) - f(\xi)\right)  g ( \xi)\big]\,. 
\en
By applying the above identity to the r.h.s. of \eqref{def_dir_form} we get
$\tilde\cE(f,g)=-\bbE[ (L f) g ]$. By item (ii) we conclude that $\tilde\cE(f,g)= - \la Lf , g\ra_{L^2(\bbP)}$, i.e. we have proved property (a).

We move to property (b), i.e.~$f\in \cD(\cE)$. Let $f_n\in \cV$ be the $n$-cutoff of $f$ given by \eqref{ncutoff}. By property (a) and Cauchy-Schwarz inequality, $\tilde\cE(f,f_n)$ is well defined and $|\tilde \cE(f,f_n)| \leq  \|  Lf \| _{L^2(\bbP)} \|  f_n \| _{L^2(\bbP)}$. Since  $\|  f_n \| _{L^2(\bbP)}\to \|  f \| _{L^2(\bbP)}<+\infty$ as $n\to+\infty$, we conclude that  $|\tilde \cE(f,f_n)|$ is bounded uniformly in $n$. Hence, to conclude that $f\in \cD(\cE)$, it is enough to show that 
\be\label{ferragosto}
\lim_{n\to+\infty}\tilde  \cE(f,f_n)=\frac{1}{2} \sum _{z\in\cN} 
	\bbE\big[ c_{0,z}(\xi) \left( f (\t_z \xi) - f(\xi)\right)^2]\,.
		\en
		To this aim we observe that $f(\xi)\geq f(\xi') \Longrightarrow f_n(\xi)\geq f_n(\xi')$ for all $\xi,\xi'\in\O$. In particular,  the functions  $h_n(\xi):=\sum_{z\in\cN}c_{0,z}(\xi)  \left( f_n (\t_z \xi) - f_n(\xi)\right)
\left( f (\t_z \xi) - f(\xi)\right)$  are nonnegative. Moreover, we observe that for $n\leq k$ we have $|f_n (\xi) - f_n(\xi')| \leq |f_k (\xi) - f_k(\xi')|$ for all $\xi,\xi'\in\O$.  This implies that  $h_n(\xi) \leq h_k(\xi)$ for all $n\leq k$.  Finally $h_n(\xi) \to h(\xi)$ for all $\xi\in \O$. We can therefore apply the monotone convergence theorem  and deduce that $\bbE[h]=\lim_{n\to+\infty }\bbE[h_n]$. Since  $\tilde\cE(f,f_n)=\bbE[h_n]/2$ while the r.h.s. of \eqref{ferragosto} equals $\bbE[h]/2$, we get \eqref{ferragosto}.

Let us now prove property (c), i.e.
\be\label{miraggio}
\cE(f,g)= - \la Lf , g\ra_{L^2(\bbP)} \qquad \forall g\in \cD(\cE)\,.
\en
Since we have shown that $f\in \cD(\cE)$, property (a) can be restated as $\cE(f,h)= - \la Lf , h\ra_{L^2(\bbP)} $  for all $ h\in \cV$.
We define $g_n\in \cV$ as the $n$--cutoff function of  $g$ (cf.~\eqref{ncutoff}).    Then, by property (a), we have $\cE(f,g_n)= - \la Lf , g_n\ra_{L^2(\bbP)}$. We now show that one can obtain \eqref{miraggio} by taking the limit $n\to +\infty$ in the above identity. To this aim recall that, when proving   that $\cV$ is a core for $\cE$, we showed that  $\cE(g-g_n,g-g_n)\to 0$. This implies that $\cE(f,g_n)\to \cE(f,g)$ (use that $\cE(f, g-g_n)\leq \cE(f,f)^{1/2} \cE(g-g_n,g-g_n)^{1/2}$), while trivially $\|g-g_n\|_{L^2(\bbP)}\to 0 $ and therefore  $\la Lf , g_n\ra_{L^2(\bbP)} \to \la Lf , g\ra_{L^2(\bbP)}$. This concludes the proof of \eqref{miraggio}.
 
We now move to step (d) and conclude. Having \eqref{miraggio} we also obtain that the map $\cD(\cE)\ni g\mapsto \cE(f,g)=- \la Lf , g\ra_{L^2(\bbP)} $ is continuous w.r.t. $\|\cdot \|_{L^2(\bbP)}$.  This observation and \eqref{miraggio} allow to apply \cite[Proposition 2.12]{R} and conclude that $f\in \cD(\bbL)$ and that $\bbL f= L f$, i.e. our criterion. We point out that \cite{R} treats also non-symmetric Dirichlet forms. On the other hand, due to \cite[Exercise 2.1]{R} and  \cite[Remark 2.2]{R}, any symmetric closed form on $L^2(\bbP)$  is a coercive closed form  on $L^2(\bbP)$ according to  \cite[Definition 2.4]{R} and therefore we can apply \cite[Proposition 2.12]{R}.

\smallskip

The last statement in Proposition \ref{zap} is an immediate consequence of the criterion we have just proved.


  %
  %
    

\begin{thebibliography}{2}


 

\bibitem{A} G.~Allaire; \emph{Homogenization and two--scale convergence}. SIAM J. Math. Anal. {\bf 23}, 1482--1518 (1992).



 
 
   \bibitem{ABSO} S. Alexander, J. Bernasconi, W. R. Schneider, R.
Orbach; {\em Excitation dynamics in random one--dimensional
systems}. Rev. Mod. Phys. {\bf 53}, 175--198 (1981)




\bibitem{AKM} S.~Armstrong, T.~Kuusi, J.-C.~Mourrat; \emph{Quantitative stochastic homogenization and large-scale regularity}. Grundlehren der mathematischen Wissenschaften Vol. {\bf 352}, Springer Verlag 2019.




\bibitem{AM}  N.\ W.\ Ashcroft, N.\ David Mermin; \emph{Solid state physics}.
 Saunders college, 1976.



 %
\bibitem{Bi} M.~Biskup; \emph{Recent progress on the random conductance model}.  
Probability Surveys, Vol.~8,  294-373 (2011).
%



 \bibitem{BP} A. Bourgeat,  A. Piatnitski; \emph{Approximations of effective coefficients in stochastic homogenization}.
  Ann. I. H. Poincar\'e  Probabilit\'es et statistiques  {\bf 40}, 153--165 (2004). 
 
 \bibitem{CI} P. Caputo, D. Ioffe;  \emph{Finite volume approximation of the effective diffusion matrix : the case of independent bond disorder}.  Ann.~Inst.~H.~Poincar\'{e} Probab.~Statist.   {\bf 39},  505--525 
 (2003).


 \bibitem{Dyre} 
J.\ C.\ Dyre, T.\ B.\ Schr{\so}der; \emph{Universality of ac conduction in disordered solids}. Rev. Mod. Phys. {\bf 72}, issue 3, 873--892 (2000).

\bibitem{F_hom} A.~Faggionato; \emph{Stochastic homogenization of random walks on  point processes}. Ann.~Inst.~H.~Poincar\'{e} Probab.~Statist. {\bf 59}, 662--705 (2023).

  \bibitem{F_resistor}  A.~Faggionato; \emph{Scaling limit of the directional conductivity of random resistor networks on simple point processes}.  Ann.~Inst.~H.~Poincar\'{e} Probab.~Statist. to appear (arXiv:2108.11258).

\bibitem{F_x_Francis} A.\ Faggionato; 
\emph{An ergodic theorem with weights and applications to random
measures, homogenization and hydrodynamics}. Stoch. Proc. their Appl. {\bf 180}, 104522 (2025).

\bibitem{FM} A.\ Faggionato, P.\ Mathieu; \emph{Linear response and Nyquist relation in periodically driven Markov processes}. In preparation. 

\bibitem{FMS} A.\ Faggionato, P.\ Mathieu, V. Silvestri; \emph{Complex mobility in periodically driven Markov processes}. In preparation.

\bibitem{FSi} A.\ Faggionato, V.\ Silvestri;
\emph{A martingale approach to time-dependent and time-periodic linear response in Markov jump processes}. ALEA, Lat. Am. J. Probab. Math. Stat. {\bf 21}, 863--906 (2024).

\bibitem{FOT} M. Fukushima, Y. Oschima, M. Takeda; \emph{Dirichlet forms and symmetric Markov processes}. Second edition. De Gruyter, Berlin, 2010.

\bibitem{GGN} N.\ Gantert, X.\ Guo, J.\ Nagel; \emph{Einstein relation and steady states for the random conductance model}. Ann. Probab. {\bf 45},  2533--2567 (2017). 

\bibitem{GMP}
N. Gantert, P. Mathieu, A. Piatnitski; \emph{Einstein relation for reversible diffusions
in a random environment}. Comm. Pure Appl. Math. {\bf 65}, 187--228 (2012).
  
\bibitem{GNO} A.\ Gloria, S.\ Neukamm, F.\ Otto; \emph{Quantification of ergodicity in stochastic homogenization: optimal bounds via spectral gap on Glauber dynamics}.
Invent. math. {\bf 199}, 455–515 (2015). 


\bibitem{JPS} R.\ Joubaud, G.A.\ Pavliotis, G.\ Stoltz; \emph{Langevin Dynamics with Space-Time Periodic Nonequilibrium Forcing}. J.~Stat.~Phys.~{\bf 158}, 1--36 (2015). 


\bibitem{KTH} R.\ Kubo, M.\ Toda, N. Hashitsume; \emph{Statistical Physics II. Nonequilibrium statistical mechanics}. Springer Verlag, Berlin, 1985.



\bibitem{MPrw} P.~Mathieu, A.~Piatnitski;  \emph{Quenched invariance principles for random walks on percolation clusters}.  Proceedings of the Royal Society A. {\bf 463} (2007) 2287--2307.


\bibitem{MP} P. Mathieu, A. Piatnitski; \emph{Steady states, fluctuation--dissipation theorems and homogenization for  reversible diffusions in a random environment}.
 Arch. Rational Mech. Anal. {\bf 230}, 277--320 (2018). 
 

\bibitem{Nu} G.~Nguetseng; \emph{A general convergence result for a functional related to the theory of homogenization}. SIAM J. Math. Anal.  {\bf 20}, 608--623 (1989).



\bibitem{O}  H. Owhadi; \emph{Approximation of the effective conductivity of ergodic media by periodization}. Probab. Theory Relat. Fields {\bf 125}, 225--258 (2003).



\bibitem{POF} M.~Pollak, M.~Ortu{\~n}o, A.~Frydman; \emph{The electron glass}.  Cambridge University Press, United Kingdom, 2013.



\bibitem{RS2}  M.~C. Reed, B. Simon, {\it Methods of modern mathematical physics. II. Fourier analysis, self-adjointness}, Academic Press, New York-London, (1975).


\bibitem{Sah} M.\ Sahimi;  \emph{Applications of percolation theory}, 2nd edition, New York,  Springer Verlag, 2023. 



\bibitem{ZP} V.V.~Zhikov, A.L.~Pyatnitskii; \emph{Homogenization of random singular structures and random measures.  (Russian) Izv. Ross. Akad. Nauk Ser. Mat. {\bf 70}, no. 1, 23--74 (2006); translation in Izv. Math. {\bf 70}, no. 1, 19--67 (2006).}

\bibitem{R} M.~R\"ockner; \emph{General Theory of Dirichlet Forms and Applications}.  In: Dell'Antonio, G., Mosco, U. (eds) Dirichlet Forms. Lecture Notes in Mathematics, vol 1563. Springer, Berlin, Heidelberg. https://doi.org/10.1007/BFb0074093


\end{thebibliography}
\end{document}